\documentstyle[12pt]{article}

\makeatletter                                            
\renewcommand{\@begintheorem}[2]{                        
\rm \trivlist \item [\hskip \labelsep {\bf #2\ \ #1.}]   
                                }                        
\makeatother                                             

\newcommand{\newsubsection}%
{{\noindent\bf\refstepcounter{subsection}\thesubsection\ \ }}

\newcommand{\newsubsubsection}%
{{\bf\refstepcounter{subsubsection}\thesubsubsection\ \ }}

\newtheorem{lem}[subsubsection]{Lemma}
\newtheorem{prop}[subsubsection]{Proposition}

\newtheorem{thm}[subsubsection]{Theorem}
\newtheorem{cor}[subsubsection]{Corollary}
\newenvironment{prf}[1]{\trivlist
\item[\hskip \labelsep{\it
#1.\hspace*{.3em}}]}{~\hspace{\fill}~$\square$\endtrivlist}

\newenvironment{proof}{\begin{prf}{Proof}}{\end{prf}}

\newcommand{\msy}{\bf}
\newcommand{\got}{{\bf}}

\newcommand{\ZZ}{{\bf Z}}
\newcommand{\DD}{{\bf D}}
\newcommand{\QQ}{{\bf Q}}
\newcommand{\AAA}{{\bf A}}
\newcommand{\GG}{{\bf G}}

\newcommand{\CC}{{\bf C}}
\newcommand{\FF}{{\bf F}}

\newcommand{\PP}{{\bf P}}

\newcommand{\cA}{{\cal A}}
\newcommand{\cC}{{\cal C}}
\newcommand{\cD}{{\cal D}}
\newcommand{\cE}{{\cal E}}
\newcommand{\cF}{{\cal F}}
\newcommand{\cG}{{\cal G}}
\newcommand{\cH}{{\cal H}}
\newcommand{\cI}{{\cal I}}
\newcommand{\cK}{{\cal K}}
\newcommand{\cL}{{\cal L}}
\newcommand{\cM}{{\cal M}}
\newcommand{\cN}{{\cal N}}
\newcommand{\cO}{{\cal O}}
\newcommand{\cP}{{\cal P}}
\newcommand{\cT}{{\cal T}}
\newcommand{\cU}{{\cal U}}
\newcommand{\cW}{{\cal W}}
\newcommand{\cX}{{\cal X}}
\newcommand{\cY}{{\cal Y}}

\def\Spec{\mathop{\rm Spec}}
\def\Diff{{\rm Diff}}
\def\Pic{{\rm Pic}}
\def\End{{{\cal E}{\it nd}}}
\def\Hom{{{\cal H}{\it om}}}
\newcommand{\dd}{{\rm d}}
\newcommand{\wP}{{\widetilde{\cP}}}
\newcommand{\bP}{{\overline{\cP}}}
\def\square{\hbox{QED}}

\begin{document}
\title{On Hitchin's Connection.}
\author{Bert van Geemen and Aise Johan de Jong}
\maketitle

{\renewcommand{\arraystretch}{1.5}

\section{Introduction}

\subsection{}
The aim of this paper is to give an explicit expression for Hitchin's 
connection in the case of rank 2 bundles with trivial determinant over
curves of genus 2.

We sketch the general situation.
Let $\pi:{\cal C}\rightarrow S$ be a family of projective smooth curves.
Let $p:\cM\rightarrow S$ be the family of moduli spaces of 
(S-equivalence classes of) rank $r$ bundles with trivial determinant
on ${\cal C}$ over $S$.
On $\cM$ we have a naturally defined determinant linebundle ${\cL}$.
The vector bundles $p_{\ast}(\cL^{\otimes k})$ on $S$ have a natural 
(projective) flat connection (such a connection also exists for other 
structure groups besides $SL(r)$).

This connection first appeared in Quantum Field Theory (conformal blocks)
and was subsequently studied by various mathematicians (see \cite{Hi} and
\cite{BM}
and the references given there). We follow the approach of Hitchin 
\cite{Hi} who showed that this connection can be constructed along the 
lines of Welters' theory of deformations \cite{W}.

In case the curves have genus 2 and the rank $r=2$,
the family of moduli spaces $p:\cM\rightarrow S$ 
is isomorphic (over $S$) to a $\PP^3$-bundle 
$$
p:\cM\cong\PP\longrightarrow S.
$$
Under this isomorphism $\cL$ corresponds to $\cO_{\PP}(1)\otimes p^*{\cal 
N}$ for some line bundle ${\cal N}$ on $S$. Giving a 
projective connection on $p_{\ast}\cL^{\otimes k}$ is then the same 
as giving a projective connection on $p_*\cO_{\PP}(k)$
$$
\nabla^{(k)}:p_\ast\cO_{\PP}(k)
\longrightarrow p_\ast\cO_{\PP}(k)\otimes\Omega^1_S.
$$

\subsection{}
In the first part of the paper (Section \ref{johan}) 
we give a characterization of
Hitchin's connection (strictly speaking, we define a new connection which
is almost characterized by Hitchin's requirements). 

We first develop a bit of general theory
concerning heat operators. A heat operator is a certain type of map
$D:p^*(\cT_S)\rightarrow \Diff^{(2)}_{\cM}(\cL^{\otimes k})$.
Thus a vector field $\theta$ 
on $S$ gives a second order differential operator $D(\theta)$ on 
$\cL^{\otimes k}$, which determines a map 
$D(\theta) : p_{\ast}\cL^{\otimes k} \to p_{\ast}\cL^{\otimes k}$. 
These then determine a connection on $p_{\ast}\cL^{\otimes k}$.
The symbol of $D$ is a map $\rho_D:p^*(\cT_S)\rightarrow S^2\cT_{\cM/S}$.

Hitchin showed that if a heat operator $D$ determines the desired 
connection on $p_{\ast}\cL^{\otimes k}$, then its symbol map:
$$
\rho_D^s:p^*(\cT_S)\longrightarrow S^2\cT_{\cM^s/S},
$$
where $\cM^s$ is the locus 
of stable bundles in $\cM$, must be equal to a map 
$1/(k+2)\rho^s_{Hitchin}$.
The map $\rho^s_{Hitchin}$ is determined by the deformation theory of 
a typical fibre $C$ and of the stable bundles on $C$.
He also proved that
$\rho^s_{Hitchin}$ extends to a 
$\rho_{Hitchin}:p^*(\cT_S)\rightarrow S^2\cT_{\cM/S}$
for bundles of rank $r$ on curves of genus $g$, 
except maybe for the case $r=2,\;g=2$ (!) (in that case the non-stable 
bundles have codimension 1 in $\cM$, whereas in all other cases they have 
codimension at least two). We show, using simple bundles,
that $\rho^s_{Hitchin}$ also extends to a $\rho_{Hitchin}$
in that case (Theorem \ref{symext}).

{}{}From a study of heat operators and their symbols we conclude that 
we can find a (projective)
heat operator $D$ with the right symbol over all of $\cM$
(Corollary \ref{exisD}). Such a $D$ need not be unique however.

The two-torsion group scheme ${\cal H}:=Pic^0_{\cC/S}[2]$ 
acts on $\cM$ (tensoring 
a rank two bundle by a line bundle of order two doesn't change the 
determinant). Then ${\cal H}$ acts (projectively)
on $\cL^{\otimes k}$, this action can be linearized to the action of
a central extension ${\cal G}$ of ${\cal H}$ on $\cL^{\otimes k}$
$$
1\longrightarrow {\bf G}_m\longrightarrow {\cal G}
\longrightarrow {\cal H}\longrightarrow 0.
$$

We show that there is a unique (projective) heat operator 
$$
\DD_k:p^*(\cT_S)\longrightarrow \Diff^{(2)}_\cM(\cL^{\otimes k})\qquad
{\rm with}\quad \rho_{\DD_k}=\mbox{${1\over {k+2}}$}\rho_{Hitchin}
$$ 
and which commutes with the action of ${\cal G}$ 
(Proposition \ref{existence+heat}).
The (projective) connection determined by $\DD_k$ is called
Hitchin's connection (cf.\ \ref{H22}).

\subsection{}
To be able to compute $\DD_k$ (and to establish 
some of the results mentioned 
above), we use a map
$f:\Pic^0_{\cC/S}\rightarrow \cM$ over $S$ whose image is 
the family of Kummer surfaces ${\cal K}\rightarrow S$.
We show that the pull-back of $\DD_k$ along $f$ determines
the connection on $q_*f^*\cL^{\otimes k}$ given by the heat 
equations on (abelian) theta functions. (The connection defined by
$\DD_k$ on $p_*\cL^{\otimes k}$ does not pull-back to the `abelian'
connection however.)

We show that $\DD_k$ is characterized by this property. Moreover, we 
express this property in terms 
of the (local) equation defining $\cK\subset \cM$ (\ref{crit} $(Eq.)$).

\subsection{}\label{defP}
To write down a completely explicit connection (on $p_*\cO_{\PP}(k)$)
we must now choose a 
certain $S$ and a family of curves $\cC$ over $S$. We take $S=\cP$,
the configuration space of $2g+2$ points in $\CC^{2g+2}$ and take the 
obvious family of hyperelliptic curves over it:
$$
\begin{array}{cll}
{\cal C}& & \supset {\cal C}_z:\quad y^2=\prod_i(x-z_i)\\
\Big\downarrow&&\phantom{\subset}\Big\downarrow \\
{\cal P}&:=\CC^{2g+2}-\{(z_1,\dots,z_{2g+2}):
\; z_i\neq z_j\;{\rm if}\;i\neq j\;\},\quad&\ni z
\end{array}
$$
We are in the case $g=2$, but the connections we define might
be interesting for any $g$.

The pull-back of the bundle $p:\PP\rightarrow \cP$ along an unramified 
$2^4:1$-Galois cover $\tilde{\cP}\rightarrow \cP$ with group $H=(\ZZ/2\ZZ)^4$
is trivial. The bundle $\PP$ can thus obtained as a quotient of
$\tilde{\cP}\times \PP^3$ by $H$ (\ref{PPdescent}).
Therefore, locally 
(in the complex or the etale topology), 
the vector bundles $p_*\cO_{\PP}(k)$ and the trivial
bundle $S_k\otimes \cO_\cP$ are canonically
isomorphic, with $S_k=H^0(\PP^3,\cO(k))$,
the vector space of homogenous polynomials of degree $k$ in 4 variables.
Thus it suffices to write down the connection on the trivial bundle which
corresponds to Hitchin's connection on $p_*\cO_{\PP}(k)$.

This connection on $S_k\otimes \cO_\cP$ is given by 
$$
\nabla:S_k\otimes\cO_\cP\longrightarrow S_k\otimes\Omega_\cP,\qquad
\nabla(w\otimes f)=w\otimes\dd f-\mbox{$1\over {k+2}$}\sum_{i\neq 
j}M_{ij}(w)\otimes{f{\dd z_i}\over{z_i-z_j}},
$$
for certain endomorphism $M_{ij}\in {\rm End}(S_k)$.
We give two descriptions of the endomorphisms $M_{ij}$.
One is based on the Lie algebra $so(2g+2)$ and its (half)-spin 
representations. We will in fact identify $S_k=S^kV(\omega_{g+1})$
where $V(\omega_{g+1})$ is a half spin representation of $so(6)$.

The other description uses the Heisenberg group action and we sketch 
that description here.
The family ${\cal C}$ has $2g+2$ sections which are the Weierstrass points 
on each fiber:
$$
P_i:\cP\longrightarrow {\cal C},\qquad z\longmapsto (x,y)=(z_i,0).
$$
Each $P_i-P_j$ is a section in $\Pic^0_{\cC/S}[2](S)\cong(\ZZ/2\ZZ)^{2g}$.
There is a central extension $G$ of this group which acts on the $S_k$'s
(the groupscheme $\cG$ is a twist of the constant group scheme $G$, these
group schemes are isomorphic over $\tilde{P}$).
Let $\pm U_{ij}\in\mbox{End}(S_1)$ be the endomorphisms induced by $P_i-P_j$ with
the property that $U_{ij}^2$ is the identity. Since 
$\mbox{End}(S_1)=S_1\otimes S_1^\ast$ and $S_1^\ast$ 
may be identified with the
derivations on $S:=\oplus S_k$ and the composition $U_{ij}U_{ij}$ may
be considered as a second order differential operator
which acts on any $S_k$. The symbol (that is, the degree two part) of 
this operator is, up to a factor $-1/16$,
the desired $M_{ij}$ (which, in this sense, does not 
depend on $k$).

The idea for trying these $\Omega_{ij}$'s comes from \cite{vGP}.
There the Hitchin map, a Hamiltonian system on the cotangent bundle
$T^*\cM^s$, is studied. 
The Hitchin map is basically the symbol of the heat operator $\cD_k$ 
(cf.\ \cite{Hi} $\S$4).
Among the results of \cite{vGP} is an explicit description of
the Hitchin map for $g=2$ curves (obtained from line geometry in $\PP^3$)
and implicitly in that description is
the one for the $M_{ij}$ given here.

\subsection{}
Now that the connection is determined it is interesting to consider its 
monondromy representation. In the case of rank two bundes, this has been 
considered by Kohno \cite{K2}. We verify that our connection gives has the
same local monodromy as Kohno's representation in Subsection 
\ref{localmon}. 
Kohno (and the physicists) use (locally on the base it 
seems) the Knizhnik-Zamolodzhikov
equations (which define a flat connection over $\cP$
using the Lie algebra $sl(2)$ \cite{Ka}, Chapter XIX).
It is not clear to us how the KZ-equations are related to our connection
(which is based on $so(2g+2)$).
We do not know which representations of the braid group one
obtains from our connections.

\section{Characterizing Hitchin's connection}\label{johan}

\subsection{Introduction}

We recall some deformation theory as described by Welters \cite{W}
and relate it to connections. 
Then we introduce heat operators  and show how they define
(projective) connections on certain bundles.
Next we give a criterion for the existence of a heat equation with a 
given symbol and in Section \ref{modbun} we recall the definition of the
symbol of Hitchin's connection. In Section \ref{g=2} we combine these 
results and give criteria (in \ref{crit} and \ref{H22}) 
which determine Hitchin's connection.

\subsection{Deformation theory}\label{def}

\subsubsection{Deformations}

Let $X$ be a smooth projective variety over $\CC$ with tangent bundle 
$T_X$.
Let $L$ be an invertible $\cO_X$-module on $X$ and let $s\in H^0(X, L)$.
We are going to classify first order deformations of $X$, of the
pair $(X,L)$ and of the triple $(X,L,s)$. See \cite{W} for proofs.

Isomorphism classes of first order deformations $X_\epsilon$ of $X$ 
(over $\Spec \CC[\epsilon]/(\epsilon^2)$) are classified by $H^1(X,T_{X})$.
We write $[X_\epsilon]$ for the class in $H^1(X,T_{X})$ of the deformation
$X_\epsilon$; similar notation will be used throughout.

The first order deformations $(X_\epsilon,L_\epsilon)$ of the pair $(X,L)$ 
are classified by $H^1(X,\Diff^{(1)}_X(L))$. Here $\Diff^{(1)}_X(L)$ is the 
sheaf of first order differential operators on $L$. This sheaf sits in an 
exact sequence:
$$
0\longrightarrow \cO_{X}\longrightarrow 
\Diff^{(1)}_X(L) \stackrel{\sigma}{\longrightarrow} T_{X}\longrightarrow 0
$$
where the map $\sigma:\Diff^{(1)}_X(L)\rightarrow T_{X}$ gives the symbol of the 
operator. This sequence gives a map $\alpha:H^1(X,\Diff^{(1)}_X(L))\to H^1(X,T_X)$,
which maps $[(X_\epsilon,L_\epsilon)]$ to $[X_\epsilon]$.

Finally, we consider first order deformations of the triple $(X,L,s)$.
The evaluation map $d^1s:\Diff^{(1)}_X(L)\rightarrow L,\;D\mapsto Ds$ gives
a complex:
$$
0\longrightarrow \Diff^{(1)}_X(L)\stackrel{d^1s}{\longrightarrow} L
\longrightarrow 0.
$$
Let ${\msy H}^1(d^1s)$ be the first hypercohomology group of this 
complex. This is the space classifying the isomorphism classes
of deformations 
$(X_\epsilon,L_\epsilon,s_\epsilon)$ of $(X,L,s)$ (\cite{W}, Prop.\ 1.2).
The spectral sequence connecting hypercohomology to cohomology gives a 
map $\beta: {\msy H}^1(d^1s) \to H^1(\Diff^{(1)}_X(L))$, which maps
$[(X_\epsilon,L_\epsilon,s_\epsilon)]$ to $[(X_\epsilon,L_\epsilon)]$.

Any element of ${\bf H}^1(d^1s)$ may be represented by a Cech cocycle.
Let us choose an affine open covering ${\cal U} : X=\bigcup U_i$ and
write $U_{ij}=U_i\cap U_j$. A Cech cocycle is given by a
pair $(\{t_i\},\,\{D_{ij}\})$ satisfying the relations:
$$
t_i-t_j=D_{ij}(s),\quad
D_{jk}-D_{ik}+D_{ij}=0\qquad \qquad
(t_i\in L(U_i),\; D_{ij}\in \Diff^{(1)}_X(L)(U_{ij})).
$$
The maps $\beta$ and $\alpha\circ\beta$ correspond to forgetting 
first $s_\epsilon$ and then $L_\epsilon$. They can be described
as follows. The cocycle $(\{t_i\},\,\{D_{ij}\})$ is mapped to
the 1-cocycle $\{D_{ij}\}$ in the sheaf $\Diff^{(1)}_X(L)$ under $\beta$. 
This is then mapped to the 1-cocycle $\{\Diff^{(1)}(D_{ij})\}$ in
the sheaf $T_X$.

Let us make explicit a deformation $(X_\epsilon, L_\epsilon, s_\epsilon)$
associated to $(\{t_i\},\,\{D_{ij}\})$. 
Write $U_i=\Spec(A_i)$ and $U_{ij}=\Spec(A_{ij})$. Let $M_i=L(U_i)$
and $M_{ij}=L(U_{ij})$. The section $s$ gives elements $s_i\in M_i$;
note that $t_i\in M_i$ as well.
Write $\theta_{ij}=\sigma(D_{ij})\in T_X(U_{ij})$.
The 1-cocycle $\{\theta_{ij}\}$ defines
a scheme $X_\epsilon$ over $\Spec(\CC[\epsilon])$ by glueing the rings
$A_i[\epsilon]$ via the isomorphisms 
$$
A_{ij}[\epsilon]\longrightarrow
A_{ij}[\epsilon],\qquad 
f+g\epsilon\longmapsto f+(\theta_{ij}(f)+g)\epsilon.
$$
Here we view $\theta_{ij}$ as a derivation. The 1-cocycle $\{D_{ij}\}$
in $\Diff^{(1)}_X(L)$ 
gives a deformation $L_\epsilon$ of $L$ by glueing the 
$A_i[\epsilon]$-modules $M_{i}[\epsilon]$ via the
$A_{ij}[\epsilon]$-module isomorphisms:
$$
\phi_{ij}:M_{ij}[\epsilon]\longrightarrow M_{ij}[\epsilon],\qquad
m_{ij}+m'_{ij}\epsilon\longmapsto 
m_{ij}+(D_{ij}(m_{ij})+m'_{ij})\epsilon.$$
Note one has $D_{ij}(fs)=fD_{ij}(s)+\theta_{ij}(f)s$ by 
definition of the symbol $\sigma$ and $\sigma(D_{ij})=\theta_{ij}$.
The cocycle $(\{t_i\},\,\{D_{ij}\})$ in ${\bf H}^1(d^1s)$ defines a 
deformation $s_\epsilon\in H^0(X_\epsilon,L_\epsilon)$. The section 
$s_\epsilon$ is given by the family of elements:
$$
s_\epsilon:=\{s_{\epsilon,i}\}\quad(\in M_i[\epsilon])
\qquad{\rm with}\quad
s_{\epsilon,i}:=s_i+t_i\epsilon.
$$ 
The
cocycle relation $t_i-t_j=D_{ij}(s)$ shows that the 
$s_{\epsilon,i}$'s glue: $\phi_{ij}(s_{\epsilon,j})=s_{\epsilon,i}$.

\subsubsection{Second order operators}\label{secondorder}

The following construction is useful to obtain elements in 
${\msy H}^1(d^1s)$ (\cite{W}, (1.9)).
Let $\Diff^{(2)}_X(L)$ be the sheaf of second order differential operators
on $L$. The symbol gives an exact sequence which is the first row 
of the complex:
$$
\begin{array}{ccccccccc}
0&\longrightarrow&\Diff^{(1)}_X(L)&\longrightarrow&\Diff^{(2)}_X(L)
&\longrightarrow&S^2T_X&\longrightarrow&0\\
&&\downarrow&&\downarrow&&\downarrow&&\\
0&\longrightarrow&L&\longrightarrow&L&\longrightarrow&0&\longrightarrow&0,
\end{array}
$$
the vertical maps are evaluation on $s$ ($D\mapsto Ds$), thus the
first column is the complex $d^1s$ considered before.
{}From the exact sequence of hypercohomology one obtains a map
$$
\delta_s:H^0(S^2T_X)\longrightarrow {\msy H}^1(d^1s).
$$

Let us give a description in terms of
cocycles of this map. Let ${\cal U}=\{U_i\}$ be an affine
open cover of $X$ as before. Let $w\in H^0(X, S^2T_X)$. We can find
operators $D^{(2)}_i\in \Diff^{(2)}_X(L)(U_i)$ which map to
$w|_{U_i}$ in $S^2T_X(U_i)$. Then
$$
\delta_s(w)=(\{D^{(2)}_i(s)\},\{D^{(2)}_i-D^{(2)}_j\})
\qquad(\in {\msy H}^1(d^1s)).
$$ 
It is easy to verify that $\delta_s(w)$ is indeed a cocycle for the
complex defined by $d^1s$. 
Note $w=\sigma(D^{(2)}_i)=\sigma(D^{(2)}_j)$ on $U_{ij}$, 
thus $D^{(2)}_i-D^{(2)}_j$ is a
first order operator.

Note that $\beta(\delta_s(w))=\{D^{(2)}_i-D^{(2)}_j\}\;(\in \Diff^{(1)}_X(L))$ 
is {\it independent}
of the choice of the section $s$. Therefore, given $w\;(\in H^0(X,S^2T_X))$
we get a deformation
$(X_\epsilon, L_\epsilon)$ such that {\it any} section $s\in H^0(X, L)$
is deformed to a section $s_\epsilon\in H^0(X_\epsilon, L_\epsilon)$.

\subsubsection{Infinitesimal connections}\label{infcon}
Let $(X_\epsilon,L_\epsilon)$ be a deformation of $(X,L)$.
We consider a differential operator $D$ of order at most two on $L_\epsilon$
which locally can be written as (cf.\ \ref{defho} below)
$$
D=\{D_i\}\quad(\in H^0(X_\epsilon,\Diff^{(2)}_{X_\epsilon}(L_\epsilon))),
\qquad
D_i=\partial_\epsilon+R_i:M_i[\epsilon]\longrightarrow M_i[\epsilon]
$$
with  $\partial_\epsilon(n+\epsilon m)=m$ for $n,\,m\in M_i$, and
$R_i$ a
$\CC[\epsilon]$-linear map of 
order at most
$2$ on $M_i$ (so $s_i\mapsto R_i(fs_i)-fR_i(s_i)$ has order at most
$k-1$ and 
a map of order zero is by definition an $A_i[\epsilon]$-linear map).
Thus $R_i$ is a second order differential operator involving only 
deriviatives in the fiber directions (and not in the base direction 
$\epsilon$).
The symbol of such a $D$ is the second order part of the $D_i$ (and thus 
of the $R_i$). The restriction of the symbol to 
$X\;(\subset X_\epsilon)$ is:
$$
\bar{\rho}(D):=\{\sigma(\bar{R}_i)\}\quad{(\in H^0(X,S^2T_{X}))},
\qquad{\rm where}\quad R_i=\bar{R}_i+\epsilon \bar{R}'_i
$$
with $\bar{R}_i,\,\bar{R}_i':M_i\rightarrow M_i[\epsilon]$ (recall $R_i$ is
$\epsilon$-linear). 

For any $s\in H^0(X,L)$ we then have a hyper cohomology class
$$
\delta_s(w):=(\{\bar{R}_i(s_i)\},\{\bar{R}_i-\bar{R}_j\})\qquad 
(\in {\msy H}^1(d^1s)).
$$

First of all we show that the deformation $(X_\epsilon,L_\epsilon)$ is the deformation
determined by $-\beta(\delta_s(\bar{\rho}(D)))$:
$$
[(X_\epsilon,L_\epsilon)]+\beta(\delta_s(\bar{\rho}(D)))=0
\qquad (\in H^1(X,\Diff^{(1)}_X(L))).
$$
This can be verified as follows. Let $s_{\epsilon,i}\in M_i[\epsilon]$ such
that $s_{\epsilon,i}=\phi_{ij}(s_{\epsilon,j})$. Then we have
$\phi_{ij}(D_j(s_{\epsilon,j})=D_i(\phi_{ij}(s_{\epsilon,j})$.
With $s_{\epsilon,j}:=s_j+t_j\epsilon,\;(s_j\,t_j\in M_j)$ this gives
(for the `constant' term):
$$
D_{ij}(s_j)+t_j+\bar{R}_j(s_j)=
t_j+\bar{R}_i(s_j),
$$
which is zero for all such local sections iff
$$
D_{ij}+(\bar{R}_j-\bar{R}_i)=0\qquad \in 
\Diff_X^{(1)}(L)(U_{ji}).
$$
Since the $\bar{R}_j$ are local lifts of $\bar{\rho}(D)$, we get that
$\beta\delta_s(\bar{\rho}(D))=\{\bar{R}_i-\bar{R}_j\}\;(\in 
\Diff_X^{(1)}(L)(U_{ij})$.

Next we consider any deformation $(X_\epsilon,L_\epsilon,s_\epsilon)$
of $(X,L,\epsilon)$.
We say that a section $s_\epsilon=\{s_i+\epsilon t_i\}
\in H^0(X_\epsilon,L_\epsilon)$ 
is flat for $D$ if
$$
D_i(s_i+\epsilon t_i)\equiv 0\;{\rm mod}\;(\epsilon),\qquad
{\rm equivalently}\quad
t_i=-\bar{R}_i(s_i),
$$
in that case we write $\nabla(D)(s_\epsilon)=0$.
Since $[(X_\epsilon,L_\epsilon,s_\epsilon)]=(\{t_i\},\{D_{ij}\})\;(\in{\bf 
H}(d^1s))$
we get:
$$
\delta_s(\bar{\rho}(D))=-[(X_\epsilon,L_\epsilon,s_\epsilon)]\quad
{\rm iff}\quad \nabla(D)(s_\epsilon)=0.
$$

\subsubsection{A remark on symbols}\label{normalization}

We use the conventions of \cite{W} concerning
differential operators and symbols. In particular
when we write $S^2 T_X$ we mean symmetric tensors
in the sheaf of vector fields on $X$. Note that if
$X=\AAA^1_\CC$ with coordinate $t$, then the symbol
of the operator $\partial^2/\partial^2t\;(\in \Diff^{(1)}_\cO)$ is the
symmetric section $2\partial/\partial t \otimes \partial/\partial t
\;(\in S^2T_X)$. By abuse of notation we will sometimes write 
$\sigma(\partial^2/\partial^2t)=2\partial^2/\partial^2t$.

\subsection{Heat operators}

\subsubsection{Notation}

Let $p: \cX\to S$ be a smooth surjective morphism of smooth
varieties over $\CC$. Let $\cL$ be an invertible $\cO_\cX$-module
over $\cX$.

We write $\cT_\cX$ (resp.\ $\cT_S$, resp.\ $\cT_{\cX/S}$) for the
sheaf of vector fields on $\cX$ (resp.\ $S$, resp.\ $\cX$ over $S$).
We denote $\Diff^{(k)}_\cX(\cL)$ the sheaf of differential operators
of order at most $k$ on $\cL$ over $\cX$. Further, 
$\Diff^{(k)}_{\cX/S}(\cL)\subset \Diff^{(k)}_\cX(\cL)$ denotes
the subsheaf of $p^{-1}(\cO_S)$-linear operators.
We remark that $\Diff^{(k)}_\cX(\cL)$ is a coherent (left) 
$\cO_\cX$-module and that $\Diff^{(k)}_{\cX/S}(\cL)$ is a
coherent submodule. For any $k$ there is a symbol map
$$ \sigma^{(k)} : \Diff^{(k)}_\cX(\cL) \longrightarrow 
S^k\cT_\cX={\rm Sym}^k_{\cO_\cX}(\cT_\cX).$$
It maps $\Diff^{(k)}_{\cX/S}(\cL)$ into $S^k\cT_{\cX/S}$.

Consider a point $x\in \cX(\CC)$, an (\'etale or analytic)
neighbourhood $U$ of $p(x)$ in $S$, and an (\'etale or analytic)
neighbourhood $V$ of $x$ in $\cX$ such that $p(V)\subset U$.
Assume we have coordinate functions $t_1,\ldots, t_r\in \cO_S(U)$ 
and functions $x_1,\ldots, x_n\in \cO_\cX(V)$ such that 
$t_1=t_1\circ p,\ldots, t_r=t_r\circ p, x_1,\ldots, x_n$ 
are coordinates on $V$. (In the \'etale case this means that 
$\Omega^1_{U}=\bigoplus \cO_U{\rm d}t_i$
and $\Omega^1_{V/U}=\bigoplus \cO_V{\rm d}x_i$.)
Then the elements $\partial/\partial x_i$ form a basis for
$\cT_{\cX/S}|_V$ over $\cO_V$. The system $(V, t_1, \ldots,
t_r, x_1,\ldots,x_n)$ will be called a {\it coordinate patch}.
We can find coordinate patches around any point $x\in \cX(\CC)$.

\subsubsection{Definition}\label{defho}

A {\it heat operator $D$ on $\cL$ over $S$} is a (left) $\cO_\cX$-module
homomorphism
$$D : p^\ast(\cT_S)\longrightarrow \Diff^{(2)}_\cX(\cL)$$
which satisfies the following property: for any coordinate patch
$(V, t_1, \ldots,t_r, x_1,\ldots,x_n)$ as above, and any
trivialization $\cL|_V\cong \cO_V$, the operator 
$D(\partial/\partial t_1)$ has the following shape:
$$
D({\partial\over\partial t_1})= f + {\partial\over\partial t_1}
+ \sum_{i=1}^n f_i {\partial\over\partial x_i} +
\sum_{1\leq i\leq j\leq n} f_{ij} {\partial^2\over\partial x_i\partial x_j}
$$
for certain $f, f_i, f_{ij}\in \cO_\cX(V)$. 
More precisely, let $$\cW_{\cX/S}(\cL)=\Diff^{(1)}_\cX(\cL)+
\Diff^{(2)}_{\cX/S}(\cL)\subset \Diff^{(2)}_\cX(\cL)$$ be 
the subsheaf of second order differential
operators on $\cL$ whose symbol lies in 
$S^2(\cT_{\cX/S})\subset S^2(\cT_\cX)$. Note that $\Diff^{(1)}_{\cX/S}(\cL)
\subset \cW_{\cX/S}(\cL)$ and that there is a canonical exact sequence
$$ 
0\longrightarrow 
\Diff^{(1)}_{\cX/S}(\cL)\longrightarrow
\cW_{\cX/S}(\cL) \longrightarrow
p^\ast(\cT_S)\oplus S^2\cT_{\cX/S}\longrightarrow 0.
$$
Our condition on $D$ above means that $D$ is a map $D: p^\ast(\cT_S)
\to \cW_{\cX/S}(\cL)$ whose composition with the map
$\cW_{\cX/S}(\cL)\to p^\ast(\cT_S)$ is the identity.
(So $D(\partial/\partial t_1)$ has no terms involving 
$\partial/\partial t_i$ except for the
term $\partial/\partial t_1$ with coefficient 1.)

Such a map $D$ is determined by the $\cO_S$-linear map
$p_\ast D: \cT_S \to p_\ast(\Diff_\cX^{(2)}(\cL))$. Note that
$\cO_S$ is a subsheaf of $p_\ast(\Diff_\cX^{(2)}(\cL))$.
Any $\cO_S$-linear map 
$\bar D : \cT_S \to p_\ast(\Diff_\cX^{(2)}(\cL))/ \cO_S$
can locally be lifted to a map into $p_\ast(\Diff_\cX^{(2)}(\cL))$.
This means that any point $s\in S(\CC)$ has a neighbourhood $U$
such that there is an $\cO_U$-linear map $D_U : \cT_S|_U\to
p_\ast(\Diff_\cX^{(2)}(\cL))|_U$ which reduces to $\bar D|_U$.

A {\it projective heat operator $\bar D$ on $\cL$ over $S$}
is an $\cO_S$-linear map
$$ 
\bar D : \cT_S \longrightarrow p_\ast\big(\Diff_\cX^{(2)}(\cL)\big)\Big/ \cO_S
$$
such that any local lifting $D_U$ as above gives rise to a
heat operator on $\cL$ over $U$.

The {\it symbol} of a heat operator $D$ is the $\cO_\cX$-linear
map 
$$ \rho_D: p^\ast(\cT_S)\longrightarrow S^2\cT_{\cX/S}$$
given by composing $D$ with the symbol $\sigma^{(2)}$. We note
that the symbol of a projective heat operator is well-defined.

\subsubsection{Heat operators and connections}

We claim that a heat operator $D$ gives rise to a connection
$$\nabla(D) : p_\ast(\cL)\longrightarrow 
p_\ast(\cL)\otimes_{\cO_S}\Omega^1_S.$$
Indeed, suppose that $s\in p_\ast(\cL)(U)=\cL(p^{-1}(U))$
and that $\theta\in \cT_S(U)$. Then $D(p^{-1}\theta)$ is
a second order differential operator on $\cL$ over $p^{-1}(U)$.
Hence $D(p^{-1}\theta)(s)\in \cL(p^{-1}(U))=p_\ast(\cL)(U)$. By our local
description of $D$ above we see that 
$D(p^{-1}\theta)(fs)=fD(p^{-1}\theta)(s)+\theta(f)s$ for every 
$f\in \cO_S(U)$.
Hence the operation $(\theta, s)\mapsto D(p^{-1}\theta)(s)$ defines
a connection $\nabla(D)$ on $p_\ast \cL$ over $S$. In the sequel we
will often write $D(\theta)$ in stead of $D(p^{-1}\theta)$.

If we have a projective heat operator $\bar D$, the connection
$\nabla(\bar D)=\nabla(D_U)$ is only defined locally, by choosing lifts 
$D_U$.
The difference of two local lifts is given by a local section 
$\eta\in \Omega^1_S(U)$, in which case the difference of the two 
connections $\nabla(D_U)$ is
multiplication by $\eta$. Thus $\bar D$ defines unambiguously
a projective connection $\bar \nabla(\bar D)$.

\subsubsection{Heat operators and deformations}\label{heat+defo}

Let $D$ be a heat operator on $\cL$ over $S$. 
Let $0\in S(\CC)$ be a point of the base. Put $X=\cX_0$ and $L=\cL_0$. 
An element $\theta\in T_0S$ will also be considered as a morphism
$\theta:\Spec(\CC[\epsilon])\to S$. By base change we get a pair
$(X_\epsilon, L_\epsilon)$ over $\Spec(\CC[\epsilon])$. Clearly,
$w_\theta:=\rho_D(\theta)$ is an element of $H^0(X, S^2T_X)$. It can be 
seen by working through the definitions that 
$[(X_\epsilon, L_\epsilon)]+\beta(\delta_s(w_\theta))=0$ for any section
$s\in H^0(X, L)$. (Note that in Subsection \ref{secondorder}
we showed that $\beta(\delta_s(w_\theta))$ was independent of $s$.)
Let $s_\epsilon\in H^0(X_\epsilon, L_\epsilon)$ lift $s$. One can show 
that $[(X_\epsilon, L_\epsilon, s_\epsilon)]=-\delta_s(w_\theta)$
if and only if the section $s_\epsilon$ is horizontal for the
connection $\nabla(D)$.

\subsubsection{Projective heat operators and change of line 
bundle}\label{heat+change}

Let $D$ be a heat operator on $\cL$ over $S$.
Let $g\in H^0(S, \cO_S^\ast)$ be an invertible function on $S$.
Multiplication by $p^{-1}(g)$ defines an invertible operator on $\cL$,
also denoted by $g$.
For any local section $\theta$ of $\cT_S$, we have the obvious relation
$$ 
g^{-1}\circ D(\theta) \circ g = g^{-1}\theta(g) + D(\theta).
$$
We conclude that $g^{-1}\circ D \circ g$ is a heat operator and that
the projective heat operators $\bar D$ and 
$\overline{g^{-1}\circ D \circ g}$ are equal. (The difference
is given by the section $g^{-1}{\rm d}g={\rm d}\log g$ of $\Omega^1_S$.)

Suppose that $\cM$ is an invertible $\cO_S$-module on $S$.
The above implies that the set of projective heat
operators on $\cL$ over $S$ can be identified canonically with the
set of projective heat operators on $\cL'=\cL\otimes_{\cO_\cX}p^\ast(\cM)$.
Indeed, choose a covering of $S$ on whose members $U_i$ the line bundle
$\cM$ becomes trivial. This identifies the heat operators on
$\cL$ and $\cL'$ over $U$. The difference in the local identifications
is given by the 1-forms ${\rm d}\log g_{ij}\in \Omega^1_{U_{ij}}$.

\subsubsection{Heat operators and flatness}\label{heat+flat}

We say that a heat operator $D$ on $\cL$ over $S$ is {\it flat} if
$D([\theta,\theta'])=[D(\theta),D(\theta')]$ for any two local 
sections $\theta, \theta'$ of $\cT_S$ on $U\subset S$. It suffices
to consider local vector fields $\theta, \theta'$ on $S$ with
$[\theta,\theta']=0$ and to check that $[D(\theta),D(\theta')]=0$.
We remark that in this case the operator $[D(\theta),D(\theta')]$
is a section of $\Diff^{(3)}_{\cX/S}(\cL)$ over $p^{-1}(U)$.
We say that $D$ is {\it projectively flat} if we have 
$D([\theta,\theta']) = h + [D(\theta),D(\theta')]$ for some function
$h=h_{\theta,\theta'}\in \cO_S(U)$ for any $\theta,\theta'\in\cT_S(U)$.
A projective heat operator $\bar D$ is called {\it projectively flat}
if any of the local lifts $D_U$ are projectively flat.

We remark that if the heat operator $D$ is flat, then the
associated connection $\nabla(D)$ on $p_\ast\cL$ is flat as well.
In particular, if $p_\ast \cL$ is a coherent $\cO_S$-module,
then this implies that $p_\ast \cL$ is locally free. In the
same vein we have that a projectively flat projective
heat operator $\bar D$ defines a projectively flat projective connection
$\bar\nabla(\bar D)$ on $p_\ast\cL$.

\subsubsection{Heat operators with given symbol}\label{heat+symbol}

Let $\cX\to S$ and $\cL$ be as above, and assume that we are given
a $\cO_\cX$-linear map 
$
\rho: p^\ast\cT_S\rightarrow  S^2\cT_{\cX/S}.
$
The question we are going to study is the following:

When can we find a (projective) heat operator $D$ on $\cL$
over $S$,
$
D:p^\ast \cT_S\rightarrow \cW_{\cX/S}\quad(\subset
\Diff^{(2)}_{\cX}(\cL))$
with $\rho_D=\rho$?

So the 
 symbol $\rho_D$ (which is the compsition of $D$ with the map
$\sigma^{(2)}:\Diff^{(2)}_{\cX}(\cL))\rightarrow S^2\cT_{\cX/S}$)
should be equal to the given $\rho$.

Let us define a number of canonical maps associated to the
situation. First we have the standard exact sequence
$$
 0 \longrightarrow \cT_{\cX/S} \longrightarrow
\cT_\cX \longrightarrow p^\ast \cT_S \longrightarrow 0.
$$
This gives us a $\cO_S$-linear map (the {\it Kodaira-Spencer map})
$$ 
\kappa_{\cX/S} \ \ :\ \  \cT_S \longrightarrow R^1p_\ast \cT_{\cX/S}.
$$
Next we have the exact sequence
$$ 
0 \longrightarrow \cT_{\cX/S}\longrightarrow
\Diff^{(2)}_{\cX/S}(\cL)\big/ \cO_\cX \longrightarrow
S^2\cT_{\cX/S} \longrightarrow 0.
$$
This gives an $\cO_S$-linear map
$$ 
\mu_\cL\ \  :\ \  p_\ast S^2\cT_{\cX/S} \longrightarrow 
R^1p_\ast \cT_{\cX/S}.
$$
The invertible $\cO_\cX$-module $\cL$ is given by an element
of $H^1(\cX,\cO_\cX^\ast)$. Using the map ${\rm d}\log : \cO^\ast_\cX
\to \Omega^1_{\cX/S}$ we get an element of $H^1(\cX, \Omega^1_{\cX/S})$,
and hence a section $[\cL]$ in $H^0(S, R^1p_\ast(\Omega^1_{\cX/S}))$.
This will be called the cohomology class of $\cL$.
By cupproduct with
$[\cL]$ we get another map $p_\ast S^2\cT_{\cX/S} \to R^1p_\ast 
\cT_{\cX/S}$,
$w\mapsto w\cup [\cL]$. (There is a natural pairing
$S^2\cT_{\cX/S}\otimes \Omega^1_{\cX/S}\to \cT_{\cX/S}$.)
It is shown in \cite[1.16]{W} that we
have $\mu_\cL(w) = - w\cup [\cL] + \mu_{\cO}(w) $.

Now let $\theta$ be a local vector field on $S$, i.e., 
$\theta\in \cT_S(U)$, with $U$ affine. 
We want to find a $D(\theta)\in \cW_{\cX/S}(\cL)$ with symbol 
$\rho(\theta)$.
Recall that $\cW_{\cX/S}(\cL)$ is an extension:
$$
0\longrightarrow 
\Diff^{(1)}_{\cX/S}(\cL)\longrightarrow
\cW_{\cX/S}(\cL) \longrightarrow
p^\ast(\cT_S)\oplus S^2\cT_{\cX/S}\longrightarrow 0.
$$
The obstruction against finding a lift
$D(\theta)\;(\in \cW_{\cX/S}(\cL))$ 
of the section $p^{-1}(\theta)\oplus\rho(p^{-1}\theta)$ 
of $p^\ast \cT_S\oplus S^2\cT_{\cX/S}$ over 
$p^{-1}(U)$
is an element $o(\theta, \cL)$ 
in $H^1(p^{-1}(U), \Diff^{(1)}_{\cX/S}(\cL))=
H^0(U, R^1p_\ast(\Diff^{(1)}_{\cX/S}(\cL)))$ 
(recall that $U$ was assumed affine). 
There is an exact sequence
$$
 0 \longrightarrow \cO_{\cX} \longrightarrow
\Diff^{(1)}_{\cX/S}(\cL)\longrightarrow \cT_{\cX/S}
\longrightarrow 0.
$$
We remark that the image of the class $o(\theta, \cL)$
in $R^1p_\ast\cT_{\cX/S}(U)$ is the section 
$\kappa_{\cX/S}(\theta) + \mu_\cL(\rho(\theta))$.
The first condition that has to be satisfied for $\bar D$
to exist is therefore: 
$$
 \kappa_{\cX/S} + \mu_\cL \circ \rho = 0. 
\leqno{(*)}
$$
If this condition is satisfied, then we can find a lift of 
$p^{-1}(\theta)\oplus\rho(p^{-1}\theta)$ to an element
$\overline{D}_\theta$ in the sheaf $\cW_{\cX/S}(\cL)/\cO_\cX$.
This element is well defined up to addition of a section
of $\cT_{\cX/S}$. Hence we get a second obstruction in
$$
 {\rm Coker}\Big(p_\ast \cT_{\cX/S} 
{\stackrel{[\cL]\cup}{\longrightarrow}}
R^1p_\ast \cO_X\Big) \leqno{(**)}
$$

Let us assume for the moment that both obstructions vanish for any
$\theta$ as above. Then we can locally find a $D(\theta)$ lifting the
element $p^{-1}(\theta)\oplus\rho(p^{-1}\theta)$. Thus, if we
have a basis $\theta_1,\ldots,\theta_r$ of $\cT_S$ over $U$,
then we can define $D_U$ by the formula
$D_U(\sum a_i \theta_i)=\sum a_i D(\theta_i)$
for certain choices of the elements $D(\theta_i)$.
If there is some way of choosing the elements $D(\theta)$ uniquely
up to an element of $p^{-1}\cO_S(U)$, e.g.\ if the map in $(**)$
is injective and $p_\ast(\cO_\cX)=\cO_S$, then $D_U$ determines
a unique projective heat operator over $U$. In this case the
map $\rho$ determines a 
unique projective heat operator on $\cL$ over $S$.

Here is another set of hypotheses which imply the existence of heat 
operators with given symbol.
Let $0\in S(\CC)$ be a point. Put $X=\cX_0$ and $L=\cL_0$
as in \ref{heat+defo}. Let $\theta\in T_0S$ and consider
$(X_\epsilon, L_\epsilon)$ over $\Spec(\CC[\epsilon])$
as in \ref{heat+defo}. We know from \ref{heat+defo}
that $[(X_\epsilon, L_\epsilon)]+\beta(\delta_s(\rho(\theta)))=0$
if a heat operator exists. Suppose that
\begin{enumerate}
\item $R^1p_\ast \cO_\cX$ and $R^1p_\ast\cT_{\cX/S}$
are locally free on $S$, and have fibres $H^1(X, \cO_X)$
and $H^1(X, T_X)$ at any point $0$ in $S(\CC)$.
\item for any $0$ and $\theta$ as above, we have
$[(X_\epsilon, L_\epsilon)]+\beta(\delta_s(\rho(\theta)))=0$.
\end{enumerate}
If these conditions are satisfied, then we can at least find
local liftings $D_U$ of $\rho$ to a heat operator. Again, we will
get a ``canonical'' projective heat operator if there is a
``preferred'' way of choosing the local lifts $D(\theta_i)$.

\subsubsection{Heat operators on abelian schemes}\label{heat+abelian}

In this subsection, we let $\cX\to S$ be an abelian scheme
with zero section $o: S\to \cX$. We assume the line bundle
$\cL$ is relatively ample on $\cX$ over $S$, i.e., $\cL$
defines a polarization of $\cX$ over $S$. We will show that
there exists a unique projective heat operator on $\cL$
over $S$, see \cite{W}.

We first remark that in this case the map $(**)$ of
\ref{heat+symbol} is an isomorphism. This implies that
$\Gamma(p^{-1}(U), \Diff^{(1)}_{\cX/S}(\cL) )=\cO_S(U)$
for any $U\subset S$. Hence, by the discussions in
\ref{heat+symbol} we get a unique heat operator as
soon as we have a (symbol) map $\rho$ satisfying $(**)$.

We have $\mu_\cO=0$ in this situation, as we can let
elements of $S^2\cT_{\cX/S}$ act by translation invariant
operators (see \cite[1.20]{W}). 
Hence, $\mu_\cL(-)= - [\cL] \cup (-) $. Now we use that
$$p_\ast (\otimes^2\cT_{\cX/S}) 
\stackrel{[\cL] \cup} {\longrightarrow}
R^1p_\ast\cT_{\cX/S}\cong o^\ast \cT_{\cX/S}\otimes R^1p_\ast\cO_\cX$$
is an isomorphism. Consider $\rho:=([\cL] \cup)^{-1} \circ \kappa_{\cX/S}$.
It is well known that this is a map into $p_\ast S^2\cT_{\cX/S}$
and it solves $(*)$ by definition; it is of course the
unique solution to $(*)$. The unique 
projective heat operator $\bar D=\bar D_{\cL}$ with $\rho=\rho_{\bar D}$
is called the projective heat operator associated to $(\cX/S,\cL)$.

We make the remark that if $\cL' = \cL^{\otimes n}$,
then we have $ \rho_{\bar D'} = ({1/n}) \rho_{\bar D}$.
It is known that these projective heat operators are
projectively flat. This can be seen in a number of ways;
one way is to show that $p_\ast(\Diff^{(3)}_{\cX/S}(\cL))=\cO_S$
in this case, compare \ref{heat+flat}.

Uniqueness of the construction of $\bar D$ implies that
$\bar D$ commutes (projectively) with the action of the theta-group
${\cal G}(\cL)$. In particular the projective monodromy
group of the (projectively flat) connection $\bar \nabla(\bar D)$
on $(p_\ast \cL)_0$ normalizes the action of ${\cal G}(L)$.

\subsubsection{Heat operators and functoriality}\label{heat+funct}

Let $\cX\to S$ and $\cL$ as above, and let $D$ be a heat operator on
$\cL$ over $S$. Note that for any morphism of smooth schemes
$\varphi:S'\to S$ there is a pullback heat operator
$D'=\varphi^\ast(D)$ on the pull back $\cL'$ on $\cX'=S'\times_S\cX$.
We leave it to the reader to give the precise definition.
This
pullback preserves flatness. If
$\varphi^\ast(p_\ast\cL)\cong p'_\ast\cL'$, then the connection
$\nabla(D')$ is the pullback of the connection $\nabla(D)$. 
In addition pullback is well defined for
projective heat operators and preserves projective flatness.

If $\psi:\cU\to \cX$ is an \'etale morphism of schemes
over $S$, then the heat operator $D$ induces a
heat operator $\psi^\ast(D)$ on $\psi^\ast(\cL)$
over $S$. 

\subsubsection{Heat operators and compatibility}\label{heat+comp}

Suppose we have a second smooth surjective morphism $q:\cY\to S$
and that we have a morphism $f: \cY \to \cX$ over $S$. We will
assume that $f$ is a submersion. This means that for any point 
$y\in \cY(\CC)$ there is a coordinate patch $(V, t_1, \ldots, t_r,
x_1,\ldots,x_c, y_1,\ldots, y_m)$ of $f(y)$ in $\cX$ such that
$(f^{-1}(V), t_1, \ldots, t_r, y_1=y_1\circ f,\ldots, y_m=y_m\circ f)$
is a coordinate patch for $y$ on $\cY$ and $x_i\circ f=0$.

Let $\cL$ be a line bundle on $\cX$ as before. 
Let us say that a (projective) heat operator $D$ on $\cL$
over $S$ is {\it (weakly) compatible} with $f:\cY\to \cX$ if
whenever we have a coordinate patch $(V, t_1, \ldots, t_r,
x_1,\ldots,x_c, y_1,\ldots, y_m)$ as above, and a trivialization
$\cL|_V\cong \cO_V$, then the operator $D(\partial/\partial t_1)$
has the following symbol:
$$ \rho_D({\partial\over\partial t_1})= \sum\nolimits_{i\leq j=1}^m f_{ij} 
{\partial^2\over\partial y_i\partial y_j} + \sum\nolimits_{i=1}^c x_i\; 
\Xi_i$$
for certain $f_{ij}\in \cO_\cX(V)$ and $\Xi_i\in S^2\cT_{\cX/S}(V)$.
More precisely, this means there exists a 
$\rho=\rho_{f,D} : q^\ast \cT_S \to S^2\cT_{\cY/S}$ 
such that the following diagram commutes:
$$\matrix{f^\ast(\rho_D) &:& f^\ast p^\ast \cT_S & \longrightarrow &
 f^\ast S^2\cT_{\cX/S}\cr
&&||&&\uparrow\cr
\rho_{f,D}& : & q^\ast \cT_S&\longrightarrow & S^2\cT_{\cY/S}\cr}$$

Suppose that $D'$ is a heat operator on $f^\ast\cL$ on $\cY$
over $S$. We say that $D'$ is {\it (weakly) compatible} with
$D$ if $D$ is (weakly) compatible with $f$ and we have
$\rho_{D'}=\rho_{f,D}$. In this case it is not true in general
that $D(\theta)(s)|_\cY=D'(\theta)(s|_\cY)$.

Assume that $D$ is compatible with $f$. Let us write
$N_\cY\cX$ for the generalized normal bundle of $\cY$
in $\cX$, i.e., $N_\cY\cX={\rm Coker}(\cT_{\cY/S}\to f^\ast\cT_{\cX/S})$.
Choose a local coordinate patch
$(f^{-1}(V), t_1, \ldots, t_r, y_1=y_1\circ f,\ldots, y_m=y_m\circ f)$
as above. Write the operator $D(\partial/\partial t_1)$ in the form 
$$D({\partial\over\partial t_1})= f + {\partial\over\partial t_1} + 
\sum\nolimits_{i=1}^c f_i {\partial\over\partial x_i} + 
\sum\nolimits_{i=1}^m g_i {\partial\over\partial y_i} +
\sum\nolimits_{i\leq j=1}^m f_{ij} {\partial^2\over\partial y_i\partial 
y_j} 
+ \sum\nolimits_{i=1}^c x_i\; \Xi_i, $$
for certain local functions $f, f_i, g_i, f_{ij}$ and second order
operators $\Xi_i$. It can be seen (and we leave this to the
reader) that the class of $\sum f_i {\partial\over\partial x_i}$
in $N_\cY\cX (f^{-1}(V))$ is independent of the choice of the
local trivialization and coordinates. The upshot of this is that
if $D$ is compatible with $\cY\to \cX$, then there is a
first order symbol
$$ \sigma_{f,D} : q^\ast\cT_S \longrightarrow N_\cY\cX.$$

We say that $D$ is {\it strictly compatible} with $f: \cY \to \cX$
if $D$ is compatible with $f$ and the symbol $\sigma_{f,D}$ is zero.
We remark that this notion is well defined for a projective heat
operator as well. 

Assume that $D$ is strictly compatible with $\cY\to\cX$.
It follows from the local description of our compatibility of
restrictions that $D$ induces a heat operator $D_\cY$
on the invertible $\cO_\cY$-module $f^\ast\cL$ over $S$.
It is characterized by the property
$D(\theta)(s)|_\cY=D_\cY(\theta)(s|_\cY)$ for any
local section $\theta$ of $\cT_S$ and any  local section
$s$ of $\cL$.
It is clear that the symbol of $D_\cY$ is equal to
$\rho_{f,D}$, hence that $D_\cY$ is compatible with $D$. 
We say that a heat operator $D'$ on $f^\ast\cL$ is {\it strictly
compatible} with $D$ if $D$ is strictly compatible with
$f:\cY\to \cX$ and $D'=D_\cY$. In this case the restriction
map
$$ p_\ast \cL \longrightarrow q_\ast f^\ast \cL$$
is horizontal for the connections induced by $D$ and $D'=D_\cY$.
Note $D_\cY$ is flat if $D$ is flat.
Similar remarks hold for projective heat operators and projective
flatness.

\subsection{Moduli of bundles}\label{modbun}

\subsubsection{The symbol $\rho_{Hitchin}$ in terms of 
moduli of bundles}\label{Hitchin-symbol}

In this subsection we explain how to get a symbol map $\rho$
as in \ref{heat+symbol} in the case of the relative moduli space
of rank 2 bundles of a family of curves.

Let $C$ be a smooth projective curve over $\Spec(\CC)$, and 
let $E$ be a stable invertible $\cO_C$-module of rank 2
with ${\rm det}(E)\cong \cO_C$. This gives a point
$[E]\in \cM_C(\CC)$. There are canonical identifications
$$
T_{[E]}\cM_C=H^1(C,\End_0(E))\qquad{\rm and}\quad
T_{[E]}^*\cM_C=H^0(C,\End_0(E)\otimes \Omega^1_C)
$$
where $\End_0(E)$ denotes the sheaf of endomorphisms of 
$E$ with trace zero. A cotangent vector 
$\phi\in T^*_{[E]}\cM_C$  
thus corresponds to a homomorphism of $\cO_C$-modules 
$\phi : E\rightarrow E\otimes \Omega^1_C$. 
(Note that this is a Higgs field on the bundle $E$.) 
Composing and taking the trace gives a symmetric bilinear pairing
$$
\matrix{T_{[E]}^*\cM_C \times T_{[E]}^*\cM_C &
\longrightarrow & H^0(C,(\Omega^1_C)^{\otimes 2})\cr
(\phi , \psi) & \longmapsto & {\rm Trace}(\phi\circ\psi).\cr} 
$$
We dualize this and use Serre duality to obtain a map
$$
\rho_{C,E}:H^1(C,T_C)\longrightarrow S^2T_{[E]}\cM_C.
$$

Let us make the observation that the construction above can be
performed in families. Let $\pi:\cC \to S$ a smooth projective family
of curves over the scheme $S$ smooth over $\Spec(\CC)$. Let
$p : \cM \to S$ be the associated family of moduli spaces of
rank 2 semi-stable bundles with
trivial determinant up to $S$-equivalence. 
We remark that $p$ is a flat projective morphism, not smooth in general. 
However, the open part of stable bundles $\cM^s\subset \cM$
is smooth over $S$. We denote this smooth morphism
by $p^s: \cM^s\to S$.

Now let $T\to S$ be a morphism of finite type. We denote by
an index ${}_T$ base change to $T$.
Let $\cE$ be a locally free sheaf of rank 2 on $\cC_T$.
For a point $0\in T(\CC)$ we put $C=\cC_0$ and $E=\cE_0=\cE|_C$. 
Suppose that for any $0$ we have (a) ${\rm det}(E)\cong \cO_C$, and
(b) the bundle $E$ is stable. In this case we get an $S$-morphism
$t: T\to \cM^s$. Analogously to the above
we have a canonical isomorphism
$$ t^\ast(\Omega^1_{\cM^s/S})\cong 
(\pi_T)_\ast(\End_0(\cE)\otimes_{\cO_{\cC_T}}\Omega^1_{\cC_T/T}). $$
Using the same pairing and dualities 
as above we get an $\cO_T$-linear map
$$ R^1\pi_\ast \cT_{\cC_T/T} \longrightarrow 
S^2\Big(t^\ast(\cT^1_{\cM^s/S})\Big).$$
If we compose this with the pullback to $T$ of the Kodaira-Spencer map 
$\cT_S\to R^1\pi_\ast \cT_{\cC/S}$ of the family $\cC\to S$, then we
get a map 
$$
 \rho_{\cE/\cC_T} : \cT_S\otimes\cO_T \longrightarrow 
S^2\Big(t^\ast(\cT^1_{\cM^s/S})\Big).
$$
If there existed a universal bundle $\cE^{univ}$ over $\cM^s$,
then we would get a ``symbol''
$$
 \rho^s_{Hitchin} = \rho_{\cE^{univ}}: (p^s)^\ast \cT_S \to 
S^2\cT_{\cM^s/S}.
$$
Although the universal bundle does not exist, the symbol
$\rho^s_{Hitchin}$ does: it is the unique $\cO_{\cM^s}$-linear
map such that, whenever $\cE/\cC_T$ is given, the pullback
of $\rho_{Hitchin}$ to $T$ agrees with $\rho^s_{\cE/\cC_T}$.
Uniqueness and existence follows from the existence of ``universal''
bundles \'etale locally over $\cM^s$ and the invariance
of $\rho_{\cE/\cC_T}$ under automorphisms of $\cE$ over $\cC_T$.

\subsubsection{Hitchin's results in genus at least 3}

Let us get back to our curve $C$ over $\Spec(\CC)$. We see
from the above that any deformation $C_\epsilon$ of the curve will give a 
tensor in $\Gamma(\cM_C^s, S^2\cT_{\cM_C})$. However, since $\cM_C$ 
has points which do not correspond to stable bundles 
(and these points are singular points of $\cM_C$ for $g>2$) 
such a deformation does not give a global tensor in 
$H^0(\cM_C,S^2T_{\cM_C})$. In \cite[Section 5]{Hi}, one can find a
sketch of a proof that the tensor extends if $g>2$ (or rank $>2$,
a case we do not even consider, although our present considerations
work in that case as well). 

In the situation $\cC\to S$ there is a canonical line bundle
$\cL$ on $\cM$; it can be defined by the formula
$$\cL = {\rm det}\; R({pr}_2)_\ast \cE^{univ},\ 
\quad pr_2 : \cC\times_S \cM \to \cM.$$
It is shown in \cite[Theorem 3.6, Section 5]{Hi}, that if $g>2$ the tensor
$2/(2k+\lambda)\rho_{Hitchin}$ is the symbol of a unique projective
flat heat operator $\bar D_{Hitchin}$ on $\cL^{\otimes k}$ over $S$.
(We remark that the extra factor 2 comes from our
way of normalizing symbols, see Subsection \ref{normalization}.)
In fact the arguments of \cite{Hi} prove the existence of this
heat operator over the moduli space of curves of genus $g>2$.

We will not use these results.

\subsection{The case of genus two curves}\label{g=2}

\subsubsection{Notation in genus 2 case}\label{notation-2}

We fix a scheme $S$ smooth over $\Spec(\CC)$ and a
smooth projective morphism $\pi:\cC\to S$, whose fibres 
$C$ are curves of genus 2. We write $0\in S(\CC)$ for
a typical point and $C=\cC_0$.

We introduce the following objects associated to the situation.
$\Pic^1=\Pic^1_{\cC/S}\to S$ denotes the Picard scheme of 
invertible $\cO_\cC$ bundles of relative degree 1 on $\cC$ 
over $S$. There is a natural morphism $\cC\to \Pic^1$ over $S$,
given by $P\in C(\CC)$ maps to $[\cO_C(P)]$ in 
$(\Pic^1)_0(\CC)=\Pic^1(C)$.
The image of this morphism is a relative divisor 
$\Theta^1=\Theta^1_{\cC/S}$
on $\Pic^1$ over $S$; of course $\Theta^1\cong \cC$ as $S$-schemes.
Let us denote $\alpha : \Pic^1\to S$ the structural morphism.
Note that 
$\alpha_\ast\cO_{\Pic^1}(2\Theta^1)$
is a locally free $\cO_S$-module of rank $4$ on $S$. We put 
$\PP=|2\Theta^1|$ equal 
to the projective space of lines in this locally free sheaf.
More precisely, we define
$$ 
p : \PP = \PP\bigg(
{\cal H}{\it om}_{\cO_S}\Big(\alpha_\ast\cO_{\Pic^1}(2\Theta^1)\;,\;
\cO_S\Big)\bigg) \longrightarrow S,
$$
see \cite[page 162]{Ha} for notation used. For a point $0\in S(\CC)$
we put $P=\PP_0$.

Let $T\to S$ be a morphism of finite type, and let
$\cE$ be a locally free sheaf of $\cO_{\cC_T}$-modules of rank 2.
For a point $0\in T(\CC)$ we put $C=\cC_0$ and $E=\cE_0=\cE|_C$. 
Suppose that for any $0$ we have 

\vspace{.5 \baselineskip}

(a) $\quad {\rm det}(E)\cong \cO_C$,

(b)$\quad$ the bundle $E$ is {\it semi}-stable.

\vspace{.5 \baselineskip}
 
(Compare with Subsection 
\ref{Hitchin-symbol}.) We recall the relative divisor
$\cD_\cE\subset \Pic^1_T$ associated to $\cE$. Set-theoretically
it has the following description:
$$ D_E := \cD_\cE \cap \Pic^1(C) = \{ [L]\in \Pic^1(C) :
\dim H^1(C, E\otimes L) \geq 1\}$$
To construct $\cD_\cE$ as a closed subscheme 
we may work \'etale locally on $T$. Hence we may assume that 
a Poincar\'e line bundle $\cL$ on $\cC_T\times_T \Pic^1_T$ exists,
i.e., which has relative degree 1 for 
$p_2:\cC_T\times_T \Pic^1_T\to \Pic^1_T$ and induces 
${\rm id}: \Pic^1_T \to \Pic^1_T$. We consider the line bundle
$$ \cN=\Big({\rm det} R(p_2)_\ast \Big( p_1^\ast(\cE) 
\otimes \cL \Big)\Big)^{\otimes -1}$$
on $\Pic^1_T$. This line bundle has a natural section (the theta
function) $\theta(\cE)$, which extends the section $1$ on the
open schematically dense subscheme $U\subset \Pic^1_T$ over
which the complex $R(p_2)_\ast \Big( p_1^\ast(\cE) 
\otimes \cL \Big)$ is trivial. We define $\cD_\cE$ to be the
zero set of the section $\theta(\cE)$; it is also the largest
closed subscheme of $\Pic^1_T$ over which the coherent
$\cO_{\Pic^1_T}$-module $R^1(p_2)_\ast \Big( p_1^\ast(\cE) 
\otimes \cL \Big)$ has rank $\geq 1$ (i.e., defined in terms of 
a fitting ideal of this sheaf). 

One can prove that $\cN\otimes \cO(-2\Theta^1_T)$ is isomorphic to 
the pullback of an invertible $\cO_T$-module on $T$.
Hence the divisor $\cD_\cE$ defines an $S$-morphism of $T$ 
into $\PP$:
$$ \varphi_{(T,\cE)} : T \longrightarrow \PP.$$
These constructions define therefore a transformation of
the stack of semi-stable rank 2 bundles with trivial determinant
on $\cC$ over $S$ towards the scheme $\PP$. It turns out that
this defines an isomorphism of the coarse moduli scheme
$$ \cM_\cC \longrightarrow \PP$$
towards $\PP$, see \cite{NR}.

A remark about the natural determinant bundle $\cL$ on $\cM_\cC$
(see \cite[page 360]{Hi}).
The relative Picard group of $\PP$ over $S$ is $\ZZ$ and
is generated by $\cO_\PP(1)$. Hence it is clear that
$\cL \cong \cO_\PP(n)\otimes p^\ast (\cN)$ for some
invertible $\cO_S$-module $\cN$ on $S$. The integer $n$
may be determined as follows. We know by \cite[page 360]{Hi}
that $\cO_P(n)^{-\lambda} \cong K_P$. The integer $\lambda=4$
in this case, hence $n=1$. We are going to define a
projective heat operator on $\cL^{\otimes k}$ over $S$, but we have seen
in Subsection \ref{heat+change} that this is the same
as defining a projective heat operator on $\cO_\PP(k)$.
Hence we will work with the line bundle $\cO_\PP(k)$ from now on.
(Note that as both $\cO_\PP(1)$ and $\cL$ are defined
on the moduli stack, they must be related by a line bundle
coming from the moduli stack $\cM_2$ of curves of genus 2;
however $\Pic(\cM_2)$ is rather small.)

Let $\Pic^0=\Pic^0_{\cC/S}\to S$ denote the Picard scheme of 
invertible $\cO_\cC$ bundles of relative degree 0 on $\cC$ 
over $S$. There is a natural $S$-morphism
$$ 
f : \Pic^0 \longrightarrow \PP 
$$
which in terms of the moduli-interpretations of both spaces
can be defined as follows:
$$ 
\Pic^0(C) \ni [L] \longmapsto [L\oplus L^{ -1}]\in \PP(\CC).
$$
It is clear that the divisor $D_{E}$ associated to
$E=L\oplus L^{ -1}$ is equal to 
$D_E= (\Theta^1_C + [L]) \cup (\Theta^1_C + [L^{-1}]) \subset \Pic^1(C)$.
This implies readily that $f^\ast \cO_\PP(1) \cong 
\cO_{\Pic^0}(2\Theta^0)$,
in fact $f^\ast\cO_\PP(1)$ may serve as the definition of
$\cO_{\Pic^0}(2\Theta^0)$ on $\Pic^0$.
This induces an isomorpism $\PP\cong |\cO_{\Pic^0}(2\Theta^0)|^*$ which
identifies $f$ with the natural map. Thus $f$ factors over
$\cK:=\Pic^0/\langle\pm1 \rangle$, the relative Kummer surface,
and
defines
a closed immersion (over $S$):
$$
\cK\hookrightarrow \PP.
$$

Any $E$ on $C$ that is semi-stable but not stable is an extension
of the form $0\to L\to E\to L^{-1}\to 0$, with $L$ of degree zero.
This bundle is $S$-equivalent to the bundle $L\oplus L^{-1}$.
Thus the open subscheme $\PP^s\subset \PP$ is equal to the
complement of the image of $f: \Pic^0 \to \PP$.

\subsubsection{Automorphisms}\label{autos}

The group scheme
$\cH=\Pic^0[2]$ of 2-torsion points of $\Pic^0$ over $S$
acts on 
the schemes $\Pic^1$, 
Indeed, if $\cE$ over $\cC_T$ is a family of locally free
sheaves as above and if $\cA$ is a line bundle on $\cC_T$
which defines a 2-torsion point of $\Pic^0$, then
$\cE\otimes \cA$ is another family of locally free sheaves
satisfying (a) and (b) of \ref{notation-2}. This defines an action
$\cH\times_S \PP\to \PP$ of $\cH$ on $\PP$, given
the identification of $\PP$ as the moduli scheme. It is easy to
see that this action is induced from the natural
action of $\cH$ on $\Pic^1$. Furthermore the morphism
$f : \Pic^0 \to \PP$ is equivariant with respect to
this action.

Note that the action just defined does not lift to an
action of $\cH$ on $\cO_\PP(1)$. Let us write $\cG$ for
theta group of the relatively ample line bundle
$\cO_{\Pic^0}(2\Theta^0)$ on $\Pic^0$ over $S$; this
group scheme gives the Heisenberg group for any point 
$0\in S(\CC)$. We remark that $\cG$ fits into the
exact sequence
$$
 1\longrightarrow \GG_{m,S} \longrightarrow
\cG \longrightarrow \cH \longrightarrow 0.
$$
There is a unique action of $\cG$ on $\cO_\PP(1)$
which lifts the action of $\cH$ on $\PP$ and
agrees with the defining action $\cG$ on $f^\ast\cO_\PP(1)$
over $\Pic^0$.

\subsubsection{Test families}

Let $\cE$ be a rank 2 free $\cO$-module on $\cC_T$ satisfying
(a) and (b) of \ref{notation-2}. We will say that the pair
$(T, \cE)$ is a {\it test family} if the following conditions 
are satisfied:
\begin{enumerate}
\item $T\to S$ is smooth,
\item all $E=\cE_0$ are simple vector bundles on $C=\cC_0$, and
\item the map $\cT_{T/S} \longrightarrow 
R^1{{\rm pr}_2}_\ast (\End_0(\cE))$ is an isomorphism.
\end{enumerate}
The last condition needs some clarification. The obstruction
for the locally free sheaf $\cE$ to have a connection on $\cC_T$
is an element in 
$$H^1(\cC\times_ST, \Omega^1_{\cC\times_ST}\otimes \End(\cE)).$$
We can use the maps 
$\Omega^1_{\cC\times_ST} \to {\rm pr}_2^\ast\Omega^1_{T/S}$
and $\End(\cE)\to \End_0(\cE)$ to project this to a section of
$\Omega^1_{T/S}\otimes R^1{{\rm pr}_2}_\ast (\End_0(\cE))$.
Whence the map of condition 3. This condition means
that the deformation of $E=\cE_0$ for any $0\in T$ is
versal in the ``vertical direction''.

The following lemma is proved in the usual manner, using deformation
theory (no obstructions !) and Artin approximation.

\begin{lem}\label{existence}For any $0\in S(\CC)$ and any simple
bundle $E$ on $C=\cC_0$ there exists a test family
$(T, \cE)$ such that $(E,C)$ occurs as one of its
fibres. \end{lem}

\begin{lem}\label{etale}Let $(T, \cE)$ be a test family. The 
morphism $\varphi_{(T,\cE)} : T \to \PP$ is
\'etale.\end{lem}

\begin{proof} Let $E$ be a simple rank 2 bundle on a smooth 
projective genus two curve $C$ over $\Spec(\CC)$ with 
${\rm det}(E)=\cO_C$. Let $\eta\in H^1(C, \End_0(E))$. For 
each invertible sheaf $L$ on $C$, with ${\rm deg}(L)=-1$ 
and given embedding $L\subset E$ consider the map
$$\End_0(E) \longrightarrow \Hom(L, E).$$
We have to show $(*)$: If for all $L\subset E$ as above $\eta$ maps to 0
in $H^1(C, \Hom(L,E))$, then $\eta=0$. 

Indeed, the condition $\eta\mapsto 0$ means that the trivial
deformation $L[\epsilon]$ of $L$ over 
$C[\epsilon]=C\times \Spec(\CC[\epsilon])$ can be embedded
into the deformation of $E$ given by $\eta$. If this holds
for all $L\subset E$, then the divisor $D_E$ does not move,
and we want this to imply that the infinitesimal deformation of
$E$ is trivial.

Let us prove $(*)$. We will do this in the case that $E$
is not stable, as we already have the result in the stable case,
see Subsection \ref{Hitchin-symbol}. Thus we may assume the
bundle $E$ is simple but not stable: $E$ is a nontrivial
extension
$$ 0\longrightarrow A \longrightarrow E \longrightarrow A^{-1}
\longrightarrow 0,$$
with ${\rm deg}_C A=0$, $A^{\otimes 2}\not\cong \cO_C$. In this
case it is easy to see that there exist three line bundles
$L_1\subset E$, $L_2\subset E$ and $L_3\subset E$ of degree $-1$ on $C$
such that the arrow towards $E$ in the following exact sequence
is surjective:
$$ 0\longrightarrow K \longrightarrow L_1\oplus L_2 \oplus L_3
\longrightarrow E \longrightarrow 0.$$
Here $K$ is just the kernel of the surjection. 
For example we can take $L_i$ of the form $A^{-1}(-P_i)$ for
$i=1,2$ and $L_3$ of the form $A(-P_3)$.
(This is the only thing we will need in the rest of the 
argument; it should be easy to establish this in the 
case of a stable bundle $E$ also.) This property then holds for
$L_i\subset E$ sufficiently general also; choose $L_i$
sufficiently general. Note that
$K\cong L_1\otimes L_2\otimes L_3$ as ${\rm det}(E)=\cO_C$.

Suppose we have a deformation $E_\epsilon$ of $E$ given
by $\eta$ as above, i.e., we have $L_i[\epsilon]\to E_\epsilon$
lifting $L\to E$. This determines an exact sequence
$$ 0\longrightarrow K_\epsilon \longrightarrow
L_1[\epsilon]\oplus L_2[\epsilon] \oplus L_3[\epsilon]
\longrightarrow E_\epsilon \longrightarrow 0.$$
But ${\rm det} E_\epsilon\cong \cO_\cC[\epsilon]$,
hence $K_\epsilon\cong 
L_1[\epsilon]\otimes L_2[\epsilon] \otimes L_3[\epsilon]
\cong (L_1\otimes L_2\otimes L_3)[\epsilon]$ is trivial
as well. Now note that by our general
position we have
$$ \dim H^0(C, \Hom(K, L_1))=\dim H^0(C, L_2\otimes L_3)=1.\leqno{(1)}$$
Similar for $L_2$ and $L_3$. Therefore, up to a unit in
$\CC[\epsilon]$, there is only one map
$K[\epsilon]\to L_i[\epsilon]$, lifting $K\to L_i$.
Thus the sequence $(1)$ above is uniquely
determined up to isomorphism, hence $E_\epsilon\cong E[\epsilon]$
is constant, i.e., $\eta=0$.
\end{proof}

\subsection{Results} This lemma completes the preparations.
In the remainder of this section we prove the desired results on the
existence and the characterization of the connection.

\begin{thm} \label{symext}
The Hitchin symbol $\rho^s_{Hitchin}:p^\ast\cT_S\rightarrow S^2\cT_{\PP^s/S}$ 
defined
in Subsection \ref{Hitchin-symbol} extends
to a symbol
$$ 
\rho_{Hitchin} : p^\ast \cT_S 
\longrightarrow S^2\cT_{\PP/S}.
$$
\end{thm}

\begin{proof} Let $(T, \cE)$ be a test family. We have the
morphism $\varphi=\varphi_{(T,\cE)} : T \to \PP$ which is \'etale
and wich induces therefore an isomorphism
$\cT_{T/S}\to\varphi^\ast \cT_{\PP/S}$. On the other
hand we have the isomorphism $\cT_{T/S} \to
R^1{{\rm pr}_2}_\ast (\End_0(\cE))$ of the test family $(T,\cE)$. 
We also have the
pairing as in Subsection \ref{Hitchin-symbol}
$${{\rm pr}_2}_\ast (\End_0(\cE)\otimes \Omega^1_{\cC_T/T})
\otimes {{\rm pr}_2}_\ast (\End_0(\cE)\otimes \Omega^1_{\cC_T/T}) 
\longrightarrow {{\rm pr}_2}_\ast ((\Omega^1_{\cC_T/T})^{\otimes 2}). $$
Note that ${{\rm pr}_2}_\ast (\End_0(\cE)\otimes \Omega^1_{\cC_T/T})$
is a locally free sheaf of $\cO_T$-modules as $\cE$ has simple fibres.
These maps and duality give us together with the Kodaira-Spencer map
a map
$$ \cT_S \otimes \cO_T \longrightarrow S^2\varphi^\ast \cT_{\PP/S}.$$
We remark that on the (nonempty) stable locus $T^s=\varphi^{-1}(\PP^s)$
this map is equal to the pullback of $\rho^s_{Hitchin}$. 

Thus, by Lemma \ref{etale}, we see that we can extend 
$\rho^s_{Hitchin}$ to any point
$x\in \PP(\CC)$ which is the image of some simple bundle. 
Note that a semi-stable bundle $E$, which is a 
nontrivial extension $0\to L\to E\to L^{-1}\to 0$
with $deg(L)=0$ is simple if and only if $L^{\otimes 2}\not\cong \cO_C$.
It follows from Lemma \ref{existence} and the above that we can extend 
$\rho^s_{Hitchin}$ to the
complement of the codimension 3 locus $f(\Pic^0[2])$,
and hence by Hartog's theorem it extends.\end{proof}

\begin{prop}\label{comp+abelian}
Let $f : \Pic^0 \to \PP$ be as in \ref{notation-2}, then $f$ is a submersion
on the locus $\cY:= \Pic^0 \setminus \Pic^0[2]$.
There is a 
commutative diagram
$$
\matrix{f^\ast(\rho_{Hitchin}) &:& f^\ast p^\ast \cT_S & 
\longrightarrow & f^\ast S^2\cT_{\PP/S}\cr
&&||&&\uparrow\cr
4\rho_{Ab}|_\cY& : & \cT_S\otimes \cO_\cY &\longrightarrow 
& S^2\cT_{\cY/S}\cr}
$$
where $\rho_{Ab}$ is the symbol of the abelian heat operator
on $\cO_{\Pic^0}(2\Theta^0)=f^\ast\cO_\PP(1)$. 
\end{prop}

\begin{proof}
Consider once again
a simple bundle $E$, which is given as an
extension $ 0\to L\to E\to L^{-1}\to 0$ with $deg(L)=0$. Clearly, the kernel
of the map
$$ 
H^1(C, \End_0(E))\longrightarrow H^1(C, \Hom(L, L^{-1}))
$$
represents deformations preserving the filtration on $E$,
i.e., tangent vectors along the Kummer surface. Dually, this
corresponds to the image of the map
$$ H^0(C, \Hom(L^{-1},L)\otimes\Omega^1_C)\longrightarrow
H^0(C, \End_0(E)\otimes\Omega^1_C).$$
We use the Trace-form to identify $\End_0(E)$ with its dual.
Note that there is a canonical exact sequence
$$ 
0\longrightarrow \cO_C \longrightarrow
\End_0(E)/L^{\otimes 2} \longrightarrow L^{\otimes -2}
\longrightarrow 0
$$
induced by the filtration $L\subset E$ on $E$.
Note that $\dim H^0(C, (\End_0(E)/L^{\otimes 2})\otimes \Omega^1_C)=2$,
since if it were $\geq 3$, then 
$\dim H^0(C, \End_0(E)\otimes \Omega^1)\geq 4$ contrary to
the assumption that $E$ is simple. Thus we get the exact sequence
$$
H^0(C,\Omega^1)\cong H^0(C, (\End_0(E)/L^{\otimes 2})\otimes \Omega^1_C)
\ll\!\!\longleftarrow H^0(C, \End_0(E)\otimes\Omega^1_C) \hookleftarrow
H^0(C, L^{\otimes 2}\otimes \Omega^1_C)\cong \CC.
$$
This is (canonically) dual to the sequence of cotangent spaces
$$ 0 \longrightarrow T_{[L]}\Pic^0 \longrightarrow 
T_{f([L])}(\PP) \longrightarrow (N_\cY \PP)_{[L]} \longrightarrow 0.$$
We want to see that our symbols in the point $[L]$ lie in the
space $S^2T_{[L]}\Pic^0$.

We have to study the map
$$ 
\matrix{ H^0(C, \End_0(E)\otimes\Omega^1_C)\times 
H^0(C, \End_0(E)\otimes\Omega^1_C) & \longrightarrow & 
H^0(C, (\Omega^1_C)^{\otimes 2})\cr
X\otimes \eta \times Y\otimes \omega &\longmapsto
& {\rm Trace}(XY)\otimes \eta\omega.\cr}
$$
But the elements of $H^0(C, L^{\otimes 2}\otimes \Omega^1_C)$,
resp.\ those of $H^0(C, \Omega^1_C)$ locally
look like 
$$\left( \matrix{ 0 & *\cr
0&0\cr}\right)\otimes \eta, \quad\hbox{ resp.\ }\quad 
\left( \matrix{ 1 & 0\cr
0&-1\cr}\right)\otimes \omega.$$
Therefore the pairing just reduces to twice the multiplication
pairing $H^0(C,\Omega^1_C) \times H^0(C,\Omega^1_C) \to 
H^0(C, (\Omega^1_C)^{\otimes 2})$.

This multiplication pairing, however, corresponds excactly
(via Kodaira-Spencer) to the symbol of the heat operator 
on the theta-divisor of the Jacobian $\Pic^0$ of $\cC$ over
$S$. Since we are looking at $2\Theta^0$, we get the
desired factor $4$, compare with Subsections 
\ref{heat+abelian} and \ref{heat+symbol}.
\end{proof}

\begin{cor}\label{exisD}
The symbol $\rho_{Hitchin}$
is invariant under
the action of the group scheme $\cH=\Pic^0[2]$ on $\PP$
and $S^2\cT_{\PP/S}$.
\end{cor}

\begin{proof}
This follows from the
observation that it is true for $\rho_{Ab|\cY}$ and thus for
$\rho_{Hitchin}|_{\cK}$,
the Hitchin symbol restricted to the Kummer surface.
Further,
one uses the remark that there are no nonvanishing elements
of $H^0(\PP, S^2\cT_{\PP/S})$ which vanish on $\cK$, that is 
$H^0(\PP,S^2\cT_{\PP/S}(-4))=0$.
\end{proof}

\begin{prop}\label{existence+heat}
There exists a unique projective heat operator $\bar D_{\lambda,k}$ 
(with $\lambda\in\CC^\ast$ be a nonzero complex number)
on $\cO_\PP(k)$ over $S$ with the following properties:
\begin{enumerate}
\item the symbol of $\bar D_{\lambda,k}$ is equal to 
      ${1\over 4\lambda}\rho_{Hitchin}$, and
\item the operator $\bar D_{\lambda,k}$ commutes (projectively) with
      the action of $\cG$ on $\cO_\PP(k)$.
\end{enumerate}
\end{prop}

\begin{proof}
This follows readily from the discussion in Subsection \ref{heat+symbol}.
Indeed, both obstructions mentioned there vanish in view of the
vanishing of $R^1p_\ast \cT_{\PP/S}$ and $R^1p_\ast\cO_\PP$. We have
a good way of choosing the elements $D(\theta)$: namely, we choose
them $\cG$-invariantly. This is possible: first choose arbitrary lifts,
and then average over liftings of a full set of sections of $\cH$
(this can be done \'etale locally over $S$). This determines the
$D(\theta)$ uniquely up to an element of $p^{-1}(\cO_S)$, as
$(p_\ast\cT_{\PP/S})^{\cH}=(0)$. Hence we get the desired
projective heat operator.
\end{proof}

It follows from Proposition \ref{comp+abelian} that the heat
operators $\bar D_{\lambda,k}$ are compatible with the morphism
$f : \cY \to \PP$. Also, for $k=1$ and $\lambda=1$, we see that 
$\bar D_{1,1}$ is compatible with the abelian heat operator on
$\cO_{\Pic^0}(2\Theta^0)$ restricted to $\cY$ (this is one of the
reasons for the factor ${1\over 4}$, but see below also).

\subsubsection{Hitchin's connection for genus 2 and rank 2}

The title of this subsection is somewhat misleading. As mentioned
before, in the paper \cite{Hi} there is no definition of a heat
operator in the case of genus 2 and rank 2. However, we propose
the following definition.

\subsubsection{Definition.}
\label{H22}
 The (projective) Hitchin connection on $p_\ast \cO_\PP(k)$
is the projective connection associated to the projective
heat operator $\DD_k:=\bar{D}_{(k+2), k}$. 

\

This definition makes sense for the following reason:
The symbol of the operator $\bar{D}_{(k+2), k}$ is equal
to ${2/(2k+4)}\rho_{Hitchin}$ which is equal to the symbol
that occurs in \cite[Theorem 3.6]{Hi} (the factor 2
comes from our way of normalizing symbols, see 
Subsection \ref{normalization}).
Thus the only assumption we needed in order to get 
$\DD_k$ was the assumption that it is compatible with the action
of $\cG$, see Subsection \ref{existence+heat}.

\subsubsection{How to compute the heat operators?}
\label{crit}

To determine $\DD_k$ we introduce heat operators $D_{\lambda,k}$
under the following assumptions:
\begin{enumerate}
\item The family of Kummer
      surfaces $\cK\subset \PP$ is given by one equation
      $F\in \Gamma(\PP, \cO_\PP(4))$.
\item We have chosen an integrable connection $\nabla_0$ on 
      $\cF:=p_\ast \cO_\PP(1)$ over $S$ which is equivariant for 
      the action of $\cG$.
\end{enumerate}
We remark that both conditions can be satisfied on the members of 
an open covering of $S$. Also the connection $\nabla_0$ in 2 is determined
uniquely up to the addition of a linear operator of the form
$\eta\cdot {\rm id}$ for some closed 1-form $\eta$ on $S$.
We remark that $\nabla_0$ induces a connection $\nabla_0$ on all the 
locally
free sheaves $S^k\cF=p_\ast \cO_\PP(k)$ over $S$, but these should not be
confused with the connections induced by the $D_{\lambda,k}$!
In addition $\nabla_0$ determines an integrable connection on
$\PP$ over $S$ and a lift of this to connections on
the sheaves $\cO_\PP(k)$ over $\PP$. In other words we
have rigidified $(\PP,\cO_\PP(k))$ over $S$.

There is a natural surjection
$$ S^2\cF\otimes_{\cO_S}S^2\cF^\ast \longrightarrow p_\ast 
S^2\cT_{\PP/S}.$$
which has a canonical splitting, given by decomposing 
$S^2\cF\otimes_{\cO_S}S^2\cF^\ast$ into irreducible ${\rm 
GL}(\cF)$-modules,
i.e., the unique ${\rm GL}(\cF)$-equivariant splitting.
Thus we identify a section $X$ of $S^2\cT_{\PP/S}$ with a section
$X$ of $S^2\cF\otimes_{\cO_S}S^2\cF^\ast$. Note that we can regard
sections of $S^2\cF\otimes_{\cO_S}S^2\cF^\ast$ as sections of
$\Diff^{(2)}_{\PP/S}\big(\cO_\PP(k)\big)$ in a natural manner,
by considering them as second order differential operators on the
affine 4-space $\underline{\Spec}(S^\ast\cF)$ invariant under the
scalar action of $\GG_{m,S}$. We write $E$ for the
Euler-vector field, i.e., the section of $\cF\otimes \cF^\ast$ that
corresponds to the identity map of $\cF$.

Assume we have a $\cG$-invariant section
$$
 X \in S^2\cF\otimes_{\cO_S}S^2\cF^\ast\otimes_{\cO_S}\Omega^1_S 
$$
with the following property
$$
 X \cdot F = \nabla_0(F) E \quad\hbox{in}\quad 
S^5(\cF)\otimes_{\cO_S}\cF^\ast\otimes_{\cO_S}\Omega^1_S,\leqno{(Eq.)}
$$
where $X\cdot F$ denotes contracting once.

If we have such an $X$ we define a heat operator
$$ 
D_{\lambda,k}:p^\ast \cT_S \longrightarrow \Diff^{(2)}_\PP(\cO_\PP(k)) 
$$
by the formula ($\theta$ a local section of $\cT_S$)
$$ D_{\lambda,k}(\theta) = \theta - {1\over\lambda} X_{\theta,k} .$$
Here $X_{\theta,k}$ means the following: first contract $X$
with $\theta$ to get a section $X_\theta$ of 
$S^2\cF\otimes_{\cO_S}S^2\cF^\ast$ and then consider this as a
second order differential operator $X_{\theta,k}$ of $\cO_\PP(k)$ on $X$
over $S$ by the remarks above. Note that the term $\theta$ acts
on $\cO_\PP(k)$ using $\nabla_0$.

\begin{lem} \label{lemchar}
Let $D_{\lambda,k}$ be the heat operator defined above.
\begin{enumerate}
\item
The heat operator $D_{\lambda,k}$ commutes with the action of $\cG$
on $\cO_{\PP}(k)$.
\item The
heat operators $D_{\lambda,k}$
are compatible with $f:\cY\to \PP$. If $\lambda= k$
then $D_{k,k}$ is strictly compatible with $f:\cY\to \PP$.
\end{enumerate}

\end{lem}

\begin{proof} 
The second statement holds since $X$ is $\cG$ invariant.
We now verify the last statement.
Note that $D_{k,k}$ is strictly
compatible with $\cK\to \PP$ if the operators $D_{k,k}(\theta)$
preserve the subsheaf
$\cI_\cK\cdot\cO_\PP(k)$ of $\cO_\PP(k)$. This can be seen
by the description of strict compatiblility in terms of
local coordinates given in Subsection \ref{heat+comp}.
We may check this condition on the affine 4-space
$\underline{\Spec}(S^\ast\cF)$. Thus let $G\in S^{k-4}(\cF)$
be a section, let $\theta$ be a local vector field
on $S$ and consider
$$ \nabla_{0,\theta}(GF) - {1\over k}X_\theta(GF)$$
We have to show that this is divisible by $F$.

We may assume that $G$ is horizontal for $\nabla_0$,
as we can find a horizontal basis for $S^{k-4}(\cF)$
locally. Thus we get
$$ G\nabla_{0,\theta}(F) - {1\over k} X_\theta (GF)=
G\nabla_{0,\theta}(F) - {1\over k} X_\theta (G) F - 
{1\over k}(X\cdot F)_\theta(G) - {1\over k} G X_\theta(F).$$
Here we have used a general formula for the application of a 
second order operator like $X_\theta$ on a product like $FG$.
Note also that $X(F)$ is simply the contration of $X\cdot F$,
hence by $(Eq.)$ we get
$$  G\nabla_{0,\theta}(F) - {1\over k} X_\theta (G) F - 
{1\over k} \nabla_{0,\theta}(F)E(G) - 
{1\over k} G E(\nabla_{0,\theta}(F)) .$$
Next, we use that acting by $E$ on a homogeneous polynomial
gives degree times the polynomial:
$$ G\nabla_{0,\theta}(F) - {1\over k} X_\theta (G) F - 
{1\over k}\nabla_{0, \theta}(F) (k-4)G - 
{1\over k} G 4 \nabla_{0,\theta}(F) = - {1\over k} X_\theta (G) F$$
which is divisible by $F$. Thus $D_{k,k}$ is indeed
stricly compatible with $f : \cY\to \PP$.
\end{proof}

\begin{thm}\label{269}
 Let $X$ be a solution to the equation
$(Eq.)$ above and let $D_{\lambda,k}$ be the heat operators defined in this 
subsection. Then
\begin{enumerate}
\item
The projective heat operator defined by $D_{\lambda,k}$
is the operator $\bar D_{\lambda,k}$ 
defined in Subsection \ref{existence+heat}. 
Thus $D_{k+2,k}$ defines the Hitchin connection on $\cO_{\PP}(k)$.
\item
The heat operators
$D_{\lambda,k}$ are projectively flat for any $k$ and $\lambda$.
\item
For $\lambda=k$, the operator $D_{k,k}$ is strictly
compatible with the abelian projective heat operator on
$\cO_{\Pic^0}(\cO(2k\Theta^0))$ over an open part of $\Pic^0$.
\end{enumerate}
\end{thm}

\begin{proof} Consider first the case of $D_{k,k}$ as defined in
\ref{crit}. We have seen above that it is strictly compatible
with $f : \cY \to \PP$. Therefore, by Subsection \ref{heat+comp}
it defines a heat operator $D_k$ on $f^\ast\cO_\PP(k)$ over $\cY$.
In view of Hartog's theorem this extends uniquely to a heat operator
$D_k$ of $\cO_{\Pic^0}(\cO(2k\Theta^0))$ on $\Pic^0$ over $S$. We have
seen in Subsection \ref{heat+abelian} that there is a unique such
heat operator, whose symbol is $1/k$ times the symbol $\rho_{Ab}$
of $D_1$. Combining this with the assertions in Subsection
\ref{comp+abelian} we see that $(1/4k)\rho_{Hitchin}-\rho_{D_{k,k}}$
vanishes along $\cK$ and hence is zero. This proves that $D_{k,k}$
agrees with the projective heat operator $\bar D_{k,k}$
defined in Subsection \ref{existence+heat}.

Now it follows from the transformation behaviour of the symbol
of $\bar D_{\lambda,k}$ and $D_{\lambda,k}$ that these agree
for arbitrary $\lambda$. Thus we get the agreement stated
in the theorem.

To see that these heat operators are projectively flat, we argue 
as follows. Let $\theta, \theta'$ be two local commuting vector 
fields on $S$. We have to see that the section
$[D_{\lambda,k}(\theta), D_{\lambda,k}(\theta')]$ of 
$\Diff^{(3)}_{\PP/S}(\cO_\PP(k))$ lies in the subsheaf
$p^{-1}(\cO_S)$. Again, for $\lambda=k$ this section
is zero when restricted to $\cK$, and again this implies
that the section is zero in that case. Now, let us compute
$$[D_{\lambda,k}(\theta), D_{\lambda,k}(\theta')]=
[\theta - {1\over\lambda} X_{\theta,k}, 
\theta' - {1\over\lambda} X_{\theta',k}]
= {1\over \lambda}\Big([\theta, X_{\theta',k}] - 
[\theta', X_{\theta,k}]\Big) +
{1\over \lambda^2} [X_{\theta,k}, X_{\theta',k}].$$
The index ${}_k$ is now superfuous, as we can see this as
an expression in $S^2\cF\otimes_{\cO_S} S^2\cF^\ast \oplus
S^3\cF\otimes_{\cO_S} S^3\cF^\ast$. We know that this expression
is zero for any $k$ with $\lambda=k$ considered as a
third order operator on $\cO_\PP(k)$. This implies that
both terms are zero, hence the expression is zero for
any $\lambda$ and any $k$. (We remark for the doubtful
reader that we will verify the flatness also by a direct computation
using the explicit description of the operator.)
\end{proof}

\section{The flat connection}\label{secheat}

\subsection{} We introduce certain families of flat connections on 
(trivial) bundles over the configuration space $\cP$ (cf.\ \ref{defP}).
These are defined by representations of the Lie algebra $so(2g+2)$.
Then we derive a convenient form for the equation of the family of Kummer 
surfaces $\cK\hookrightarrow\PP$ (in fact for any $g$ we point out a 
specific section $P_z$ of a trivial bundle over $\cP$). In Theorem 
\ref{thmcon}
we show that this equation is flat for one of our connections. This will
be important in identifying Hitchin's connection in the next section.

\subsection{Orthogonal groups}\label{ortgrp}
\subsubsection{The Lie algebra}
Let $Q=x_1^2+\ldots+x_{2g+2}^2$, then the (complex)
Lie algebra $so(Q)=so(2g+2)$ is:
$$
so(2g+2)=\{A\in End(\CC^{2g+2}):\; {}^tA+A=0\;\},\qquad{\rm let}\quad
F_{ij}:=2(E_{ij}-E_{ji})\;(\in\;so(2g+2)),
$$
where $E_{ij}$ is the matrix whose only non-zero
entry is a $1$ in the $(i,j)$-th position (the commutator of such
matrices is 
$[E_{ij},E_{kl}]=\delta_{jk}E_{il}-\delta_{li}E_{kj}$). 
The alternating matrices $F_{ij}$ 
 (with
$1\leq i<j\leq 2g+2$) are a basis of $so(2g+2)$. Note $F_{ij}=-F_{ji}$.
  These matrices satisfy the relations:
$$
[F_{ij},\,F_{kl}]=
\left\{\begin{array}{rcl}
0&{\rm if}&i,\,j,\,k,\,l\quad\mbox{are distinct} ,\\
2F_{il}&{\rm if}& j=k.
\end{array}\right. 
$$

\subsubsection{The universal envelopping algebra} 
Recall that for a Lie algebra ${\got g}$ 
the tensor algebra $T({\got g})$ and the universal enveloping algebra 
$U({\got g})$
are defined by:
$$
T({\got g}):=\CC\oplus{\got g}\oplus {\got g}\otimes{\got g}\ldots,
\qquad
U({\got g}):=T({\got g})/I
$$
where $I$ is the ideal generated by all elements of the form
$x\otimes y-y\otimes x-[x,y]$ with $x,\,y\in {\got g}$.
Lie algebra representations $\rho:{\got g}\rightarrow End(V)$ 
(so $\rho([x,y])=\rho(x)\rho(y)-\rho(y)\rho(x)$ correspond
to representations $\tilde{\rho}:U({\got g})\rightarrow 
End(V)$ of associative algebras with $\rho(1)=id_V$. Given $\rho$ one 
defines 
$\tilde{\rho}(x_{i_1}\otimes x_{i_2}\ldots\otimes x_{i_k}):=
\rho(x_{i_1})\rho(x_{i_2})\ldots \rho(x_{i_k})$ which is well-defined 
because we work modulo $I$.

\subsubsection{Definitions}\label{oneform}
We define elements in $U(so(2g+2))$ by:
$$
\Omega_{ij}:=F_{ij}\otimes F_{ij}\quad{\rm mod}\;I\qquad
(\in U(so(2g+2))),\qquad{\rm note}\quad \Omega_{ij}=\Omega_{ji}.
$$
In particular, for any Lie algebra representation $\rho: 
so(2g+2)\rightarrow
W$ we now have endomorphisms 
$\tilde{\rho}(\Omega_{ij})=\rho(F_{ij})^2:W\rightarrow W$
which satisfy the same commutation relations as the $\Omega_{ij}$.

For $\lambda\in\CC^*$ we define an $U(so(2g+2))$-valued one form 
$\omega_\lambda$ 
 on $\cP$:
$$
\omega_\lambda:=
\lambda^{-1} \left(\sum_{i\neq j}
{{{\Omega_{ij}} {\rm d}z_i}\over{z_i-z_j}}\right);\qquad
{\rm let}\quad \tilde{\rho}(\omega_\lambda):=
\lambda^{-1} \left(\sum_{i\neq j}{{\tilde{\rho}(\Omega_{ij})\dd z_i} \over
{z_i-z_j}}\right)
$$
where $\tilde{\rho}:U(so(2g+2))\rightarrow End(W)$ is a representation.
Note $\tilde{\rho}(\omega_\lambda)\in \mbox{End}(W)\otimes\Omega^1_\cP$.
 
The one form
$\omega_\lambda$
 defines a connection on the trivial bundle $W\otimes_\CC\cO_\cP$ by:
$$
\nabla^W_\lambda:W\otimes\cO_\cP\longrightarrow W\otimes \Omega^1_\cP,
\qquad fw\longmapsto w\otimes{\rm d}f-f\tilde{\rho}(\omega_\lambda)(w)
\qquad\qquad(f\in\cO_\cP,\;w\in W).
$$
The covariant derivative of $f\otimes w$ with respect to the vector field 
$\partial/\partial z_i$ on $\cP$ is the composition of $\nabla^W_\lambda$
with contraction:
$$
\left(\nabla^W_\lambda(f\otimes w)\right)_{\partial/\partial z_i}=
w\otimes {\partial f\over\partial z_i} -\sum_{j,\,j\neq 
i}\tilde{\rho}(\Omega_{ij})(w)\otimes{1\over z_i-z_j}\qquad\qquad
(\in W\otimes\cO_\cP).
$$

\begin{prop} \label{defcon}
Let $\rho:so(2g+2)\longrightarrow End(W)$ be a Lie algebra representation, 
then $\nabla^W_\lambda$ is a flat connection
on $W\otimes\cO_\cP$ for any $\lambda\in\CC^*$. 
\end{prop}

\begin{proof}
Obviously $\nabla^W_\lambda$ is a connection.
It is well known (and easy to verify) that 
the connection defined by $\omega$ is flat (so 
${\rm d}\omega+\omega\wedge\omega=0$) if the following 
{\it infinitesimal 
pure braid relations}  in $U(so(2g+2))$ are satisfied:
$$
[\Omega_{ij},\Omega_{kl}]=0,\qquad [\Omega_{ik},\Omega_{ij}+\Omega_{jk}]=0,
$$
where $i,\,j,\,k,\,l$ are distinct indices (cf.\ \cite{Ka}, Section XIX.2).

To check these relations we use that $X\otimes Y=Y\otimes X+[X,Y]$ in 
$U(so)$.
The first relation is then obviously satisfied. We spell out the second.
Consider first:
$$
\begin{array}{rcl}
F_{ik}\otimes F_{ij}\otimes F_{ij}&=&
(F_{ij}\otimes F_{ik}+[F_{ik},F_{ij}])\otimes F_{ij}\\
&=&F_{ij}\otimes F_{ik}\otimes F_{ij}-F_{kj}\otimes F_{ij}\\
&=&F_{ij}\otimes(F_{ij}\otimes F_{ik}-F_{kj})+F_{jk}\otimes F_{ij}\\
&=&F_{ij}\otimes F_{ij}\otimes F_{ik}+
(F_{ij}\otimes F_{jk}+F_{jk}\otimes F_{ij}),
\end{array}
$$
so we have $[F_{ik},F_{ij}\otimes F_{ij}]=D_1$ with $D_1$
symmetric in the indices $i$ and $k$. Similarly we get:
$$
\begin{array}{rcl}
F_{ik}\otimes( F_{ik}\otimes F_{ij}\otimes F_{ij})&=&
F_{ik}\otimes( F_{ij}\otimes F_{ij}\otimes F_{ik}+D_1)\\
&=&(F_{ik}\otimes F_{ij}\otimes F_{ij})\otimes F_{ik}
+F_{ik}\otimes D_1\\
&=& F_{ij}\otimes F_{ij}\otimes F_{ik}\otimes F_{ik}+(D_1\otimes F_{ik}+
F_{ik}\otimes D_1).
\end{array}
$$
Thus $[\Omega_{ik},\Omega_{ij}]=D_1\otimes F_{ik}+
F_{ik}\otimes D_1$ which is antisymmeric in $i$ and $k$. Therefore
$$
[\Omega_{ik},\Omega_{ij}]=-[\Omega_{ki},\Omega_{kj}]=
-[\Omega_{ik},\Omega_{jk}]
$$
which proves the second infinitesimal braid relation.
\end{proof}

\subsubsection{Differential operators}\label{symcon}
Given a Lie algebra representation $\rho:{\got g}\rightarrow {\rm End(S_1)}$
where $S_k\subset S:=\CC[\ldots,X_i,\ldots]$ is the subspace of homogeneous 
polynomials of degree $k$, there is a convenient way to determine the 
representation $\rho^{(k)}:g\rightarrow {\rm End}(S_k)$ induced by $\rho$.
If $\rho(A)=(a_{ij})$ w.r.t.\ the basis $X_i$ of $S_1$, then we
define
$$
L_A:=\sum_{i,j}a_{ij}X_i\partial_j\qquad(\in {\rm End}(S)),\qquad
{\rm with}\quad \partial_j(P):={\partial P\over\partial X_j}
$$
for $P\in S$. Then obviously $\rho(A)(P)=L_A(P)$ if $P$ is linear and the
Leibnitz rule shows that $\rho^{(k)}(A)=L_A:S_k\rightarrow S_k$ for all 
$k$.

The composition (in ${\rm End}(S)$) of two operators is given by:
$$
(X_i\partial_j)\circ (X_k\partial_l)=
X_iX_k\partial_j\partial_l+\delta_{jk}X_i\partial_l,\qquad
{\rm thus}\quad
L_A^2=(1/2)\sigma(L_A^2)+L_{A^2},
$$
with symbol 
$\sigma(L_A^2)=2\sum_{i,j,k,l}a_{ij}a_{kl}X_iX_k\partial_j\partial_l$
(see our convention \ref{normalization}).

Assume that we have $A_{ij}\in{\got g}$ such that 
$\Omega_{ij}=A_{ij}\otimes A_{ij}$ satisfy the infinitesimal braid 
relations (cf.\ proof of the previous proposition) and that in the 
representation $\rho$ we have $\rho(A_{ij})^2=\mu I$ with $\mu\in\CC$.
Then we conclude that the operators 
$$
\sigma(L_{A_{ij}}^2):S\longrightarrow S
$$
also satisfy the infinitesimal braid relations.

\subsection{The Kummer equation.}\label{kummer}

\subsubsection{Introduction}
We recall how one can obtain the equation of the Kummer 
variety for a genus two curve, and more generally, how a hyperelliptic
curve ${\cal C}_z$ determines in a natural way a quartic polynomial 
$P_z\in S^4V(\omega_{g+1})$ where 
$V(\omega_{g+1})$ is a half spin representation of $so(2g+2)$.

The polynomial $P_z$ lies in the subrepresentation of $V(4\omega_{g+1})$ 
which also 
occurs in $S^2(\wedge^{g+1}\CC^{2g+2})$ where $\CC^{2g+2}$ is the standard
representation of $so(2g+2)$. We will exploit this fact to verify that
$P_z$ is a flat section for one of the connections introduced in 
Subsection \ref{ortgrp}.

\subsubsection{The orthogonal Grassmanian}
We define two quadratic forms on $\CC^{2g+2}$ by:
$$
Q:\quad x_1^2+\ldots+x_{2g+2}^2,\qquad
Q_z:\quad z_1x_1^2+\dots+z_{2g+2}x_{2g+2}^2,
$$
and we use the same symbols to denote the corresponding quadrics in 
$\PP^{2g+1}$. The quadric $Q$ has two rulings (families of linear $\PP^g$'s
lying on it). 
Each of these is parametrized by the orthogonal Grassmanian (spinor 
variety)
denoted by $Gr_{SO}$. 
The variety $Gr_{SO}$ smooth, projective of dimension 
$\mbox{${1\over 2}$}g(g+1)$. 

Let $V(\omega_{g+1})$ be a half spin representation of $so(2g+2)$.
There is an embedding 
$$
\phi:Gr_{SO}\longrightarrow \PP V(\omega_{g+1})^*,\qquad {\rm let}\quad
\cO_{Gr_{SO}}(1):=\phi^*\cO_{\PP V}(1).
$$
in fact, $\phi(Gr_{SO})$ is
the orbit of the highest weight vector.
The map $\phi$ is equivariant for the action of the spin group 
$\widetilde{SO}$ 
(which 
acts through a half-spin representation on $V(\omega_{g+1})$).
It induces isomorphisms:
$$
H^0(Gr_{SO},\cO(n))=V(n\omega_{g+1}).
$$

In case $g=2$, one has isomorphisms: 
$$
\phi:Gr_{SO}\stackrel{\cong}{\longrightarrow} \PP^3,\qquad
S_k=V(k\omega_{g+1})\qquad
{\rm and}\quad \widetilde{SO}_6\cong SL_4(\CC)
$$ 
and the halfspin representation is identified  
with the standard representation $SL_4(\CC)$ on $\CC^4$ (or its dual).

\subsubsection{The Pl\"ucker map}
Let $Gr_{SL}$ be the Grassmanian of $g$-dimensional subspaces of 
$\PP^{2g+1}$, we denote the Pl\"ucker map by:
$$
p:Gr_{SL}\longrightarrow \PP \wedge^{g+1}\CC^{2g+2},\qquad
\langle v_1,\ldots,v_{g+1}\rangle \longmapsto v_1\wedge\ldots\wedge 
v_{g+1}=\sum p_{i_1\ldots i_{g+1}}e_{i_1}\wedge \ldots \wedge e_{i_{g+1}}.
$$
The map $p$ is equivariant for the $SL(2g+2)$-action on both sides
and induces isomorphisms of $sl(2g+2)$-representations
$$
H^0(Gr_{SL},\cO_{Gr_{SL}}(n))=V(n\lambda_{g+1}),\qquad
{\rm with}\quad 
V(\lambda_{g+1})\cong \wedge^{g+1}\CC^{2g+2}.
$$

As a representation of $so(2g+2)$ the space $\wedge^{g+1}\CC^{2g+2}$ is
reducible, one has (cf.\ \cite{vG2} (3.7), \cite{FH}, $\S$19.2, Th.\ 19.2 
and Remarks (ii)), with $V(\omega_g$ and $V(\omega_{g+1})$ the half spin 
representations:
$$
\wedge^{g+1}\CC^{2g+2}=V(2\omega_{g})\oplus V(2\omega_{g+1}).
$$
In particular, the irreducible component $V(4\omega_{g+1})$ of
$S^4V(\omega_{g+1})$ is also a component of $S^2(\wedge^{g+1}\CC^{2g+2})$
(viewed as $so(2g+2)$ representation). Elements of $
S^2(\wedge^{g+1}\CC^{2g+2})$ can be seen as quadratic forms in the 
Pl\"ucker coordinates, when restricted to $Gr_{SO}$ they can be viewed as 
(restrictions to $Gr_{SO}$ of) quartic polynomials on $\PP 
V(\omega_{g+1})^*$.

\subsubsection{Notation}\label{bsnot}
The following lemma gives this decomposition of $\wedge^{g+1}\CC^{2g+2}$
explicitly. Let $\{e_i\}$ be the standard basis of $\CC^{2g+2}$. For any 
$$
S\subset B:=\{1,\,2,\ldots,2g+2\},\qquad |S|=g+1,
$$
we write $S=\{i_1,\,i_2,\ldots,i_{g+1}\}$ with $i_1<i_2<\ldots<i_{g+1}$
and define
$$
e_S:=e_{i_1}\wedge\ldots\wedge e_{i_{g+1}}\qquad\in \wedge^{g+1}\CC^{2g+2}.
$$
These elements give a basis of $\wedge^{g+1}\CC^{2g+2}$. 
For such an $S$ we let $S':=B-S$ be the complement of $S$ in $B$.
Writing $S'=\{j_1,\,j_2,\ldots,j_{g+1}\}$ with $j_1<j_2<\ldots<j_{g+1}$,
we define an element $\sigma_S$ in the symmetric group $S_{2g+2}$ by:
$$
\sigma_S(k):=i_k,\quad \sigma_S(g+1+k):=j_k.
$$

\begin{lem}\label{resp}
The two non-trivial $so(2g+2)$-invariant subspaces in 
$\wedge^{g+1}\CC^{2g+2}$
are:
$$
\langle\ldots, e_S+sgn(\sigma_S)i^{g+1}e_S',\ldots\rangle_{S\ni 1}\qquad
{\rm and}\qquad 
\langle\ldots, e_S-sgn(\sigma_S)i^{g+1}e_S',\ldots\rangle_{S\ni 1},
$$
where $S$ runs over the subsets of $B$ with $g+1$ elements with $1\in S$.
\end{lem}

\begin{proof} (Cf.\ \cite{FH}, $\S$19.2 Remarks(iii).)
The quadratic form $Q$ defines an $so$-equivariant isomorphism 
$$
B:\wedge^{g+1} \CC^{2g+2}\longrightarrow (\wedge^{g+1} \CC^{2g+2})^*=
\wedge^{g+1} (\CC^{2g+2})^*,\qquad
e_S\longmapsto
\epsilon_S:=\epsilon_{i_1}\wedge\ldots\wedge \epsilon_{i_{g+1}}
$$
with $\{\ldots,\epsilon_j,\ldots\}$ the basis dual to the $e_i$.
There is also a canonical (in particular  $so$-equivariant) isomorphism:
$$
C:\wedge^{g+1} \CC^{2g+2}\longrightarrow (\wedge^{g+1} \CC^{2g+2})^*,\qquad
\alpha\longmapsto [\beta\mapsto c_{\alpha,\beta}]
$$
with
$  c_{\alpha,\beta}\in\CC$ defined by
$$
\alpha\wedge\beta=c_{\alpha,\beta}e_1\wedge e_2\wedge\ldots e_{2g+2}.
$$
Thus we have an isomorphism 
$A:=B^{-1}\circ C:\wedge^{g+1} \CC^{2g+2}\rightarrow\wedge^{g+1} 
\CC^{2g+2}$ whose 
eigenspaces are $so$-invariant. These are the subspaces we have to 
determine.  

Note $e_S\wedge e_T\neq 0$ iff $T=S'$, the 
complement of $S$. One easily verifies 
that if $1\in S$ we have:
$$
e_S\wedge e_{S'}=
sgn(\sigma_S)e_1\wedge e_2\wedge\ldots e_{2g+2},\qquad
{\rm so}\quad C(e_S)=sgn(\sigma_S)\epsilon_{S'}.
$$
As $e_{S'}\wedge e_S=(-1)^{g+1}e_S\wedge e_{S'}$, we get 
$C(e_{S'})=sgn(\sigma_S)(-1)^{g+1}\epsilon_S$ (with $1\in S$).
Therefore 
$$
A(e_S)=sgn(\sigma_S)e_{S'},\qquad 
A(e_{S'})=sgn(\sigma_S)(-1)^{g+1}e_S\qquad
{\rm if}\quad 1\in S.
$$
The eigenvalues of  $A$ are thus $\pm i^{g+1}$
and the eigenvectors are $e_S\pm sgn(\sigma_S)i^{g+1}e_S'$.
\end{proof}

\subsubsection{A quartic polynomial}\label{defpz}
We consider the $g$-dimensional subspaces of $\PP^{2g+1}$ which are tangent
to $Q_z$:
$$
\bar{B}_z:=\left\{\PP^g\subset\PP^{2g+1}:\;rank(Q_z\; \mbox{restricted 
to}\;\PP^{g})\leq g\;\right\}\qquad(\subset Gr_{SL}).
$$
The subvariety $\bar{B}_z$ is defined by a quadratic polynomial $P_z$
in the 
Pl\"ucker coordinates, that is by a section in 
$H^0(Gr_{SL},\cO_{Gr_{SL}}(2))$
(\cite{vG2}, Thm 3):
$$
\bar{B}_z=Z(P_z),\qquad P_z\in 
H^0(Gr_{SL},\cO_{Gr_{SL}}(2))\subset S^2(\wedge^{g+1}\CC^{2g+2}).
$$
This section is in fact the unique $SO_z$-invariant
in $H^0(Gr_{SL},\cO_{Gr_{SL}}(2))$
and will be determined explicitly in the lemma below.
In case $g=2$ the quartic surface in $\PP V$ defined by $P_z$
is the Kummer surface of the curve $\cC_z$ (\cite{DR}).

\begin{lem}
The variety $\bar{B}_z$ is defined by 
$$
P_z=\sum_S z_Se_S^2,
$$
where $S$ runs over the subsets with $g+1$ elements of 
$B=\{1,\ldots,2g+2\}$
and
$$
z_S:=z_{i_1}z_{i_2}\ldots z_{i_{g+1}},\qquad 
e_S=e_{i_1}\wedge\ldots\wedge e_{i_{g+1}}\qquad{\rm if}\quad
S=\{i_1,\ldots,i_{g+1}\}.
$$
\end{lem}

\begin{proof}
{}From \cite{vG2}, Theorem $3'$ and its proof we know that the 
trivial representation of $so_z$ has multiplicity one in 
$H^0(Gr_{SL},\cO_{Gr_{SL}}(2))=V(2\lambda_{g+1})$.
This $sl_{2g+2}$ representation corresponds to the partition $\mu$:
$2g+2=2+2+\ldots+2$ and is thus realized as
$(\CC^{2g+2})^{\otimes 2g+2}c_\mu$ where $c_\mu$ is the Young symmetrizer,
one has in fact $V(2\lambda_{g+1})\subset S^2(\wedge^{g+1}\CC^{2g+2})$
(see \cite{FH}, $\S$15.5, p.\ 233-236). We use the standard tableau as
in \cite{vG2} (and \cite{FH}) to construct $c_\mu$.

The trivial $so_z$ sub-representation is obtained as follows.
The quadratic form $Q_z$ defines the $so_z$-invariant tensor 
$\sum_{i=1}^{2g+2}
z_ie_i\otimes e_i\in \CC^{2g+2}\otimes \CC^{2g+2}$, and we consider the
$g+1$-fold tensor product of this tensor, where we insert the $k$-th 
factor 
in the positions $k$ and $g+1+k$:
$$
\tau_z:=\sum_{i_1,\ldots,i_{g+1}=1}^{2g+2} 
z_{i_1}z_{i_2}\ldots z_{i_{g+1}}e_{i_1}\otimes e_{i_2}\otimes\ldots
\otimes e_{i_{g+1}}\otimes e_{i_1}\otimes\ldots e_{i_{g+1}}\quad
(\in (\CC^{2g+2})^{\otimes 2g+2}).
$$
Thus $\tau_z$ is $so_z$ invariant and it has the advantage that 
$\tau_zc_\mu$ is easily determined:
$$
\tau_zc_\mu=(cst)\sum_S z_S\,e_S\cdot e_S, 
$$
with $(cst)$ a non-zero integer and
where we sum over all subsets $S$ with $g+1$ elements of 
$\{1,\ldots,2g+2\}$
(In fact, $c_\mu=a_\mu b_\mu$ 
and $a_\mu$ symmetrizes the indices $1,\,g+2$; $2,\,g+3$;
$\ldots$; $g+1,\,2g+2$ of a tensor (but $\tau_z$ is already symmetric in 
these indices) and next $b_\mu$ antisymmetrizes  the first $g+1$ and last
$g+1$ indices, giving the result above).
\end{proof}

\subsection{Verifying the contraction criterium}
\subsubsection{}\label{actf}
We need to determine the images of the $\Omega_{ij}$'s in the 
representation
$S^2(\wedge^{g+1}\CC^{2g+2})$.
In the standard representation $\CC^{2g+2}$ of $so(2g+2)$,
with basis $e_i$ and dual basis $\epsilon_i$, we have:
$$
F_{ij}=2(e_i\otimes\epsilon_j-e_j\otimes\epsilon_i)\qquad
\in
(\CC^{2g+2})\otimes(\CC^{2g+2})^*=End(\CC^{2g+2}).
$$
On $\wedge^{g+1}\CC^{2g+2}$, with basis 
$e_S=e_{i_1}\wedge\ldots\wedge e_{i_{g+1}}$ as before, $F_{ij}$ acts as
$$
F_{ij}(e_{i_1}\wedge\ldots\wedge e_{i_{g+1}})=
F_{ij}(e_{i_1})\wedge e_{i_2}\wedge\ldots\wedge e_{i_{g+1}}+\ldots+
e_{i_1}\wedge\ldots\wedge e_{i_g}\wedge F_{ij}(e_{i_{g+1}}).
$$
Define a subset of $B=\{1,\ldots,2g+2\}$ by
$$
S+ij=(S\cup\{i,j\})-(S\cap\{i,j\}).
$$
Then 
$F_{ij}(e_S)$ is, up to sign, $2e_{S+{ij}}$ if $|S\cap\{i,j\}|=1$
and is zero otherwise:
$$
F_{ij}=2\sum_{S,\;|S\cap\{i,j\}|=1} t_{S,ij}e_{S+ij}\otimes\epsilon_S
$$
with $t_{S,ij}= \pm 1$.
Note that if $|S\cap\{i,j\}|=1$ then $F_{ij}^2e_S=-4e_S$ so 
$t_{S,ij}=-t_{S+ij,ij}$. 

\begin{lem}\label{cont2}
Let $W$ be a complex vector space, let $D\in W\otimes W^*$ be a linear
differential operator and let $P\in S^kW$. Then the contraction of
$P$ with the symbol of $D^2$ is given by:
$$
\sigma(D^2)\cdot P=2D(P)D\qquad \in S^{k-1}W\otimes W^*.
$$
\end{lem}

\begin{proof}
We view elements of $W$ as linear forms in variables $X_a$. Then elements
of $W^*$ are linear operators with constant coefficients and 
$D=\sum_af_a\partial_a$ with $f_a\in W$. The symbol of $D^2$ is then:
$$
\sigma(D^2)=\sum_{a,b} f_af_b\partial_a\partial_b\qquad\in S^2W\otimes 
S^2W^*.
$$

Recall that we use the convention (see \ref{normalization}):
$$
\partial_a\partial_b:=
\partial_a\otimes\partial_b+\partial_b\otimes\partial_a
$$
for elements in $S^2W^*$. The contraction of a polynomial $P\in S^kW$
with $\partial_a\partial_b\in S^2W^*$ is thus:
$$
(\partial_a\partial_b)\cdot P=
(\partial_bP)\partial_a+(\partial_aP)\partial_b\qquad \in S^{k-1}W\otimes 
W^*.
$$
Therefore
$$
\begin{array}{rcl} \sigma(D^2)\cdot P&=&
\sum_{a,b}f_af_b\left((\partial_bP)\partial_a+(\partial_aP)\partial_b\right)\\

&=&2\sum_af_a(\sum_bf_b(\partial_bP))\partial_a\\
&=&2D(P)D.
\end{array}
$$
\end{proof}

\begin{thm}\label{thmcon}
The composition, still denoted by $P_z$:
$$
\cP\longrightarrow S^2\wedge^{g+1}\CC^{2g+2}
\longrightarrow S^4V(\omega_{g+1}),\qquad 
z\longmapsto P_z,
$$
with $P_z=\sum_Sz_Se_S^2$ (as in Subsection \ref{defpz}), satisfies the 
differential equations:
$$
\left(\frac{1}{16}
\sum_{j\neq i}{\sigma\tilde{\rho}(\Omega_{ij})\over{z_i-z_j}}\right)
\cdot 
P_z+(\partial_{z_i}P_z)E=0\qquad\qquad(1\leq i\leq 2g+2)
$$
with $E=\sum e_S\otimes\epsilon_S$ the Euler vector field and 
$\tilde{\rho}:U(so(2g+2))\rightarrow {\rm End}(S^2(\wedge^{g+1}\CC^{2g+2}))$.

\noindent
Thus $P_z$ is a horizontal section for the connection $\nabla^W_\lambda$
with $W=S^4V(\omega_{g+1})$ and $\lambda=-16$.
\end{thm}

\begin{proof}
For symmetry reasons it suffices to verify the equation for $i=1$.
By Lemma \ref{cont2}
$$
\sigma\tilde{\rho}(\Omega_{1j})\cdot e_S^2=
2\rho(F_{1j})(e_S^2)\rho(F_{1j})=4e_S\rho(F_{1j})(e_S)\rho(F_{1j}).
$$
To simplify notation we write $\Omega_{ij}$ for $\sigma\tilde{\rho}(\Omega_{ij})$
and $F_{ij}$ for $\rho(F_{ij})$ (as in \ref{actf}).
Then $F_{1j}(e_S)=0$ unless $|S\cap\{1,j\}|=1$.
Assume $1\in S$ and $j\not \in S$ (so $j\in B-S$) and consider 
 the contraction of $\Omega_{1j}$ with the term 
$z_Se_S^2+z_{S+ij}e_{S+ij}^2$ from $P_z$:
$$
\begin{array}{rcl}
\Omega_{1j}\cdot (z_Se_S^2+z_{S+1j}e_{S+1j}^2)&=&
8(t_{S,1j}z_S+t_{S+1j,1j}z_{S+1j})e_Se_{S+1j}F_{1j}\\
&=&8t_{S,1j}(z_1z_{\bar S}-z_jz_{\bar{S}})e_Se_{S+1j}F_{1j}\\
&=&(z_1-z_j)\Omega_{1j}\cdot z_{\bar{S}} e_S^2,
\end{array}
$$
where ${\bar S}:=S-\{1\}$. Therefore we get:
$$
{\Omega_{1j}\over{z_1-z_j}}\cdot P_z=
\Omega_{1j}\cdot\left(\sum_{S\ni 1,\,j\not\in S} z_{\bar{S}} e_S^2\right).
$$
Summing this result over all $j$ and changing the order of summation gives:
$$
\begin{array}{rcl}
\left(\sum_{j\neq 1}{\Omega_{1j}\over{z_1-z_j}}\right)\cdot P_z&=&
\sum_{S\ni 1} z_{\bar{S}}\left(\sum_{j\in B-S} \Omega_{1j}\cdot 
e_S^2\right)\\
&=&
8\sum_{S\ni 1} z_{\bar{S}}e_S\left(\sum_{j\in B-S} t_{S,1j} 
e_{S+1j}F_{1j}\right)
\end{array}
$$
where, as before, $B=\{1,\ldots,2g+2\}$.
On the other hand, since $P_z=\sum z_Se_S^2$, we have:
$$
 \partial_{z_1}P_z=\sum_{S\ni 1}z_{\bar{S}}e_S^2.
$$
Thus the theorem follows if we prove, for all $S$ with $1\in S$:
$$
2e_SE+\sum_{j\in B-S} t_{S,1j} e_{S+{1j}}F_{1j}=0.
$$

With the definition of $F_{1j}$ we find:
$$
\begin{array}{rcl}
\sum_{j\in B-S} t_{S,1j} e_{S+1j}F_{1j}&=&
\sum_{j\in B-S} t_{S,1j} e_{S+1j}\left(\sum_{T,\;|T\cap\{1,j\}|=1} 
2t_{T,1j}e_{T+1j}\otimes\epsilon_T\right)\\
&=&
2\sum_{T}\left(\sum_{j\in B-S,\;|T\cap\{1,j\}|=1} t_{S,1j}t_{T,1j} 
e_{S+1j}e_{T+1j}\right)\otimes \epsilon_T.
\end{array}
$$
Note that $e_SE=\sum_T e_Se_T\otimes \epsilon_T$ (with $T\in B,\;|T|=g+1)$.
Comparing coefficients of $\epsilon_T$,
it remains to prove that (for all $S\ni 1$ and all $T$):
$$
e_Se_T+\sum_{j\in B-S,\;|T\cap\{1,j\}|=1} t_{S,1j}t_{T,1j} 
e_{S+1j}e_{T+1j}=0.
$$

We show that the relations for $1\in T$ follow from those with $1\not\in 
T$.
In fact, we only want this relation in $V(4\omega_{g+1})$,
so we restrict ourselves to the subspace $V(2\omega_{g+1})\subset 
\wedge^{g+1}\CC^{2g+2}$. The kernel of the restriction map is
$\langle\ldots,e_T-sgn(e_T)i^{g+1}e_{T'},\ldots\rangle$ (or the other space
in Lemma \ref{resp}; the argument we give leads to the same conlusion in 
both cases).

Consider a $T$ with $1\in T$ and $j\not\in T$. Then $1\not\in T+{1j}$ and
after restriction: 
$$
e_T=sgn(\sigma_T)i^{g+1} e_{T'},\qquad 
e_{T+1j}=sgn(\sigma_{(T+1j)'})i^{-(g+1)} e_{(T+1j)'}
$$
where $'$ stands for the complement in $B$. Since $(T+1j)'=T'+1j$, 
substituting these relations and multiplying throughout by $i^{g+1}$ we 
get:
$$
(-1)^{g+1}sgn(\sigma_T)e_Se_{T'}+
\sum_{j\in B-S,\;|T\cap\{1,j\}|=1} t_{S,1j}t_{T,1j}sgn(\sigma_{T'+1j}) 
e_{S+1j}e_{T'+1j}
$$
Next we observe that 
$$
sgn(\sigma_{T'+1j})=(-1)^{g+1}t_{T,1j}t_{T',1j}sgn(\sigma_T),
$$ 
so it suffices to consider the relations with $1\not\in T$. 
Let $T:=\{1=i_1,\ldots,i_{g+1}\},\;T':=\{j_1,\ldots,j_{g+1}\}$ with $j=j_k$
and $i_l<j<j_{l+1}$.
To get $\sigma_{T'+1j}$, 
first apply the permutation $(1\;g+2)\ldots(g+1\;2g+2)$ to $B$, the sign 
of this permutation is $(-1)^{g+1}$. Then apply $\sigma_T$. Next apply the 
cyclic
permutation $(j_1\;\ldots j_{g+1}\;1)$ (with sign $(-1)^{g+1}$) and 
finally apply the inverse of $(j_k\;\ldots j_{g+1}\;i_2\;\ldots i_l)$ 
(with sign
$(-1)^{g-k+l}$). The resulting permutation is $\sigma_{T'+1j}$ and has sign
$(-1)^{g+k+l}sgn(\sigma_T)$. On the other hand,
since $1\in T,\;j\not\in T$ we get
$F_{1j}(e_T)=-2(-1)^{l-1}e_{T+{1j}}$, so $t_{T,1j}=(-1)^{l}$ and 
$F_{1j}(e_{T'})=+2(-1)^{k-1}e_{T'+1j}$, so $t_{T',1j}=(-1)^{k-1}$. This 
gives
the formula for $sgn(\sigma_{T'+1j})$.

>{}From now on we consider only $T$'s with $1\not\in T$ and $S$ with $1\in 
S$. 
We show that the desired relation:
$$
e_Se_T+\sum_{i\in (B-S)\cap T} t_{S,1i}t_{T,1i} e_{S+1i}e_{T+1i}.
$$
 is a Pl\"ucker relation. Thus it holds in
$H^0(Gr_{SL},\cO(2))$ and therefore also upon restriction to 
$Gr_{SO}\subset Gr_{SL}$.

Let $S=\{1=i_1,\ldots,i_{g+1}\}$ and
consider
the Zariski open subset $U\subset Gr_{SL}$ of $g+1$-dimensional
subspaces $W\subset\CC^{2g+2}$ with 
Pl\"ucker coordinate $p_{S}\neq 0$. Any such $W$
has a (unique) basis $\{\ldots,w_i,\ldots\}$ with: 
$$
W=\langle w_1,\ldots,w_{g+1}\rangle,\qquad
{\rm and}\quad
(w_k)_{i_l}=\delta_{kl},\qquad 
(S=\{1=i_1,\ldots,i_{g+1}\},\;i_1<\ldots<i_{g+1})
$$
where $\delta_{kl}$ is Kronecker's delta.
Let $M$ be the $(g+1)\times(2g+2)$ matrix whose rows are the $w_i$.
Let $M_i$ be the $i$-th column of $M_W$ (note $M_{i_j}=f_j$, 
the $j$-th standard basis vector of $\CC^{g+1}$).
Then the Pl\"ucker coordinate $p_{k_1,\ldots,k_{g+1}}$ is the determinant 
of the $(g+1)\times (g+1)$ submatrix of $M$ whose $j$-th column is 
$M_{k_j}$,
we write:
$$
p_{k_1,\ldots,k_{g+1}}=det(M_{k_1},\ldots,M_{k_{g+1}})\qquad
{\rm with}\quad k_1<\ldots<k_{g+1}.
$$
 
Let $j\not\in S$, with $i_l<j<i_{l+1}$. Then $p_{S+1j}=(-1)^{l-1}(M_j)_1$
since
$$
det(f_{i_2},\!\ldots,f_{i_l},M_j,f_{i_{l+1}},\!\ldots,f_{g+1})
\!=(-1)^{l-1}det(M_j,f_{i_2},\!\ldots,f_{i_l},f_{i_{l+1}},\!\ldots,f_{g+1})
=\!(-1)^{l-1}(M_j)_1.
$$
As earlier, since $1\in S$ we have $t_{S,1j}=-(-1)^{l-1}$, hence 
$p_{S+1j}=-t_{S,1j}(M_j)_1$.

Let $T=\{j_1,\ldots,j_{g+1}\}$, let $T\cap S'=\{a_1,\ldots,a_q\}$
and let $j=j_n=a_k$.
Then $p_T=det(M_{j_1},\ldots,M_{j_{g+1}})$ is, up to sign, 
the determinant of a $q\times q$
submatrix $N$ of the $(g+1)\times q$ matrix with columns
$M_{a_1},\ldots,M_{a_q}$, we write $p_T=t\cdot det(N)$ with $t=\pm 1$.
Then 
$$
\begin{array}{rcl}
p_{T+1j}&=&det(f_1,M_{j_1},\ldots,\widehat{M_{j_n}},\ldots,M_{j_{g+1}})\\
&=&
(-1)^{j_n-1}det(M_{j_1},\ldots,M_{j_{n-1}},f_1,M_{j_{n+1}},\ldots,M_{j_{g+1}})\\
&=&(-1)^{j_n-1}(-1)^{k+1}t\cdot det(N^{1k}),
\end{array}
$$
where $N^{1k}$ is the $(q-1)\times(q-1)$ submatrix of $N$ obtained by 
deleting the first row and $k$-th column. 
Since $1\not\in T$, we have $t_{T,1j}=(-1)^{j_n-1}$ and thus
$p_{T+1j}=t(-1)^{k+1}det(N^{1k})$. 

Substituting these expressions for $p_{S+1j}$ and $p_{T+1j}$ we get:
$$
1\cdot t\cdot det(N)+\sum_{j\in \{a_1,\ldots,a_q\}}
t(-1)^k(M_j)_1det(N^{1k})
$$
which is zero in virue of a well known formula for the determinant.
\end{proof}

\section{The Heisenberg group and the Spin representation} \label{Hspin}

\subsection{}
We recall the basic facts on the Heisenberg group and we discuss the 
projective representation of its automorphism group.
In Subsection \ref{spin}  
 we relate the Heisenberg group (in its 
irreducible $2^g$-dimensional representation) with the (half) spin 
representation of the orthogonal group.
Combining this with previous results, we can finally write down 
Hitchin's connection in Subsection \ref{hitcon}.

\subsection{The Heisenberg group}
\label{Hgroup}

\subsubsection{Definitions}\label{Hdefs}
We introduce a variant of the Heisenberg group (cf.\ Subsection 
\ref{autos}).
For any positive integer $g$ we define 
a (finite) Heisenberg group $G$ by:
$$
G=G_g:=\{(t,x):\;t\in\CC,\;t^4=1,\quad
x=(\xi,\xi')\in \FF_2^g\times \FF_2^g\}
$$
with identity element $(1,0)$ and multiplication law:
$$
(t,(\xi,\xi'))(s,(\eta,\eta')):=(ts(-1)^{\xi\eta'},\xi+\eta,\xi'+\eta')
$$
with $\xi\eta':=\sum_{i=1}^g\xi_i\eta'_i$. The group $G$ has $2^{2g+2}$
elements, it is non-abelian, in fact its center is $\{(t,0)\}$.
The inverse of $(t,x)$ can be found as follows:
$$
(t, (\xi,\xi'))(s, (\xi,\xi'))=((-1)^{\xi\xi'}ts,(0,0)),\qquad
{\rm thus}\quad (t, x)^{-1}=((-1)^{\xi\xi'}t^{-1},x),
$$
So for any non-zero $x$ there is are elements $(t,x)\in G$
of order 2 (and also of order 4).

We define a symplectic form on $\FF_2^{2g}$ by:
$$
E(x,y)=\xi\eta'+\eta\xi',\qquad{\rm with}\quad
x=(\xi,\xi'),\;y=(\eta,\eta')
$$
thus $E(x,x)=0$ for all $x$ but $E$ is non-degenerate.
Note that $E$ is related to the commutator in $G$:
$$
(t,x)(s,y)(t,x)^{-1}(s,y)^{-1}=((-1)^{E(x,y)},0).
$$

\subsubsection{Representations.}\label{repH}
Let $V$ be the $2^g$ dimensional vector space of complex valued functions
$\FF^g_2\rightarrow \CC$. It has a standard basis consisting of 
$\delta$-functions 
$$
X_\sigma:\FF^g_2\longrightarrow \CC,\qquad X_\sigma(\sigma)=1,\quad
X_\sigma(\rho)=0\quad{\rm if}\;\sigma\neq\rho.
$$
The Heisenberg group $G$ has a representation $U$ on 
$V$, 
the Schr\"odinger
representation:
$$
(U(t,(\xi,\xi'))f)(\sigma):=t(-1)^{\sigma\xi'}f(\sigma+\xi),\qquad
{\rm thus}\quad 
U(t,(\xi,\xi'))X_\sigma=t(-1)^{(\sigma+\xi)\xi'}X_{\sigma+\xi}.
$$
For every $x\in\FF_2^{2g}-\{0\}$ choose a $(t_x,x)\in G$ of order 
two. We define:
$$
U_x:=U(t_x,x)\qquad (\in GL(V)),\qquad{\rm so}\quad U_x^2=I
$$
and define $U_0=I$. Then $Im(U)=\{tU_x:\;t^4=1,\;x\in\FF^{2g}_2\}$.

\subsubsection{Automorphisms.}
We define a subgroup of $Aut(G)$ by:
$$
A(G):=\left\{\phi\in 
Aut(G):\;\phi((t,0))=(t,0)\quad\forall t\right\}.
$$
The elements of $A(G)$ are the automorphisms which are the
identity when restricted to the center of $G$. For $\phi\in 
A(G)$ and $(t,x)\in G$ we can then write:
$$
\phi(t,x):=(f_\phi(x)t,M_\phi(x)),\qquad{\rm with} \quad
M_\phi:\FF_2^{2g}\longrightarrow \FF_2^{2g},\quad
f_\phi:\FF_2^{2g}\longrightarrow \CC^*
$$
(note $\phi(t,x)=\phi(t,0)\phi(1,x)=(t,0)\phi(1,x)$).
The map $M:A(G)\rightarrow Aut(\FF_2^{2g}),\;\phi\mapsto M_\phi$ 
is a homomorphism.

Assume that $\phi\in\ker(M)$, thus $M_\phi$ is
the identity.
Then one verifies that  $f_\phi$ is a homomorphism so 
we must have $f_\phi(x)=(-1)^{E(x,y)}$ for some $y=y_\phi\in\FF_2^{2g}$.
But then $\phi$ is an interior automorphism since also
$$
Int_y:(t,x)\longmapsto (1,y)(t,x)(1,y)^{-1}=((-1)^{E(x,y)}t,x).
$$
Thus $\ker(M)\cong G/Center(G)\cong\FF^{2g}_2$.

Since automorphisms preserve commutators in $G$, 
the image of $M$ lies in
$Sp(2g,\FF_2)=Sp(\FF_2^{2g},E)$. There is an exact sequence:
$$
0\longrightarrow \FF^{2g}_2\longrightarrow A(G)
\stackrel{M}{\longrightarrow}Sp(2g,\FF_2)\longrightarrow 0.
$$
(See Theorem \ref{twine} below for the surjectivity of $M$.)

\subsubsection{A projective representation.}\label{projrep}
The Schr\"odinger representation $U$
of $G$ on $V$ is the unique
irreducible representation of $G$ in which $(t,0)$ acts by
mutiplication by $t$. Given $\phi\in A(G)$, the representation
$U\circ\phi$ enjoys the same property. Hence by Schur's lemma we get
a linear map, unique up to scalar multiple,
$$
\tilde{T}_\phi:V\longrightarrow V,\qquad{\rm with}\quad
\tilde{T}_\phi U(h)=U(\phi(h))\tilde{T}_\phi,
$$
for all $h\in G$.
In this way we get a projective representation of $A(G)$:
$$
\tilde{T}:A(G)\longrightarrow PGL(V),\qquad
\phi\longmapsto \tilde{T}_\phi.
$$
Note that we may take $\tilde{T}_\phi=U_y$ when $\phi=Int_y$. 

On the other hand, as $U:G\rightarrow GL(V)$ is injective
we have $Im(U)\cong G$. Any  $T\in GL(V)$ which normalizes $Im(U)$ 
thus
defines an automorphism $\phi_T$ of $G$ (if $TU(h)T^{-1}=U(h')$ 
then $\phi_T(h):=h'$). Since the center of $G$ acts by scalar 
multiples of the 
identity, we have $\phi_T\in A(G)$. Thus we get an exact sequence
$$
0\longrightarrow \CC^*\longrightarrow Normalizer_{GL(V)}(Im(U))
\longrightarrow A(G)\longrightarrow 0.
$$

\subsubsection{Definition}
For $x=(\xi,\xi')\in \FF_2^{2g}-\{0\}$, we define the transvection
$$
T_x:\FF_2^{2g}\longrightarrow \FF_2^{2g},\qquad
y\longmapsto y+E(y,x)x.
$$
Then $T_x\in Sp(2g,\FF_2)$ and 
the (finite) symplectic group $Sp(2g,\FF_2)$ is generated
by transvections. Note that the transvections are involutions: $T^2_x=1$.

\begin{thm}\label{twine}
For $x\in \FF^2_{2g}-\{0\}$ let 
$$
\tilde{T}_x:=U_x+iI,\qquad{\rm with}\quad i^2=-1.
$$
Then $\tilde{T}_x\in A(G)$ and 
$$
M(\tilde{T}_x)=T_x\qquad\mbox{ that is:}\qquad
\tilde{T}_xU_y\tilde{T}_x^{-1}=t_{x,y}U_{T_x(y)},
$$
for all $y\in\FF_2^{2g}$ and some $t_{x,y}\in \CC$ with $t^4_{x,y}=1$.
Therefore the homomorphism $M:A(G)\rightarrow Sp(2g,\FF_2)$ is
surjective.

Moreover, in $A(G)$ we have
$\tilde{T}_x^2=Int_x$.
\end{thm}

\begin{proof} 
Since $U^2_x=I$, the eigenvalues of $U_x$ are $\pm 1$. Thus $U_x+iI$ is 
invertible (its inverse is $(1/2)(U_x-iI)$). In case $E(x,y)=0$, $U_x$ and
$U_y$ commute and $T_x(y)=y$, thus the relation holds. In case $E(x,y)=1$,
$U_xU_y=-U_yU_x$ and $T_x(y)=T_{x+y}$, and the relation holds because:
$$
\begin{array}{rcl}
\tilde{T}_xU_y\tilde{T}_x^{-1}&=&
(1/2)\left(U_xU_yU_x-i(U_xU_y-U_yU_x)+U_y\right)\\
&=&
(1/2)(-U_y-2iU_xU_y+U_y)\\
&=& \pm U_{x+y}\qquad(\in Im(U)),
\end{array}
$$
in fact, $U_{x+y}$ and $-iU_xU_y\in Im(U)$ differ by a scalar multiple $t$ 
and since both elements have order two, $t=\pm1$. 
Thus $\tilde{T}_x$ normalizes $Im(U)$
and defines an element of $A(G)$ indicated by the same symbol.
The relation shows that $M(\tilde{T}_x)=T_x$. 
Since the $T_x$ generate $Sp(2g,\FF_2)$ the map
$M$ is surjective. 

Finally $(U_x+iI)^2=U_x^2+2iU_x-I=2iU_x$, and conjugation by $U_x$ induces 
$Int_x$, which proves the last statement.
\end{proof}

\subsubsection{Example}\label{extt}
A particular case
is when $x=(0,\xi')$ with $\xi'=(1,0,\ldots,0)\in\FF_2^g$. Then on the 
basis
$X_{(0,\ldots,0)},\ldots,X_{(0,\tau)},\ldots,X_{(1,0,\ldots,0)},\ldots,
X_{(1,\tau)},\ldots$ with $\tau\in\FF^{g-1}_2$ 
we have:
$$
U_x=\pmatrix{I&0\cr 0&-I},\qquad
\tilde{T}_x':=\mbox{$1\over{1+i}$}\tilde{T}_x=\pmatrix{I&0\cr 0&iI},
$$
note that $\tilde{T}_x'$ and $\tilde{T}_x$ define the same element of 
$A(G)$.

\subsection{Notation}

\subsubsection{}\label{defB}
In dealing with hyperelliptic curves and the half-spin representation of
the orthogonal group, the following (classical) notation for points in
$\FF_2^{2g}$ is convenient (\cite{DO}, VIII.3;
\cite{M}).
Let
$$
B:=\{1,\,2,\ldots,2g+2\},\qquad{\rm then}\quad
F_B:=\left\{f:B\rightarrow \FF_2:\;\sum_{b\in B}f(b)=0\right\}\Big/
\{0,[b\mapsto 1]_{b\in B}\}.
$$
is an $\FF_2$-vector space of dimension $2g$.
For a subset $T\subset B$ with an even number elements we denote
by $x_T\in F_B$ (or simply $T$)
the element defined by the function $f$ with $f(b)=1$
iff $b\in T$. Note that $x_T=x_{T'}$ when $T'$ is the complement of 
$T$ in $B$, moreover
$$
x_T+x_S=x_R\qquad{\rm with}\quad R=T+ S:=(T\cup S)-(T\cap S).
$$

\subsubsection{}\label{notiso}
We fix the following isomorphism of $\FF_2$-vector
spaces, and identify them in this way in the remainder
of the paper
(\cite{DO}, VIII.3, Lemma 2, but note we interchanged 
$(e_i,0)\leftrightarrow (0,e_i)$):
$$
F_B\stackrel{\cong}{\longrightarrow} \FF_2^{g}\times\FF^{g}_2,\qquad
x_{\{2i-1,\,2i\}}=(0,e_i),\quad
x_{\{2i,\,2i+1,\ldots,2g+1\}}=(e_i,0)\qquad (1\leq i\leq g).
$$
The symplectic form $E$ on $\FF_2^{g}\times\FF^{g}_2$  
can be now be easily computed on $F_B$ by:
$$
E(x_T,x_S)=|T\cap S| \;\;{\rm mod}\,2,\qquad
{\rm thus}\quad
E(x_{ij},x_{kl})=
\left\{\begin{array}{rcl}
0&{\rm if}&i,\,j,\,k,\,l\quad\mbox{are distinct} ,\\
1&{\rm if}& i<j=k<l
\end{array}\right.
$$
where we write $x_{ij}$ for $x_{\{i,j\}}$ etc. For $g=1$ we have
$$
x_{12}=x_{34}=(1,0),
\quad x_{13}=x_{24}=(1,1),\quad
x_{14}=x_{23}=(0,1).
$$
For $g=2$ one has 15 non-zero points $x_{ij}=x_{klmn}$ when 
$\{i,j,k,l,m,n\}=B$:
$$
x_{12}=x_{3456}=
((0,0),(1,0)),\quad
x_{2345}=x_{16}=
((1,0),(0,0)),\quad
x_{26}=x_{1345}=
((1,0),(1,0)).
$$

\subsection{The Spin representation}\label{spin}

\subsubsection{}
The two half spin representations
of $so(2g+2)$ are each realized on a $2^g$-dimensional vector space.
We recall, using the Clifford algebra, how they can be constructed
using the Heisenberg group. 
We will have to consider the Heisenberg group 
$G_{g+1}$
 which acts on a $2^{g+1}$-dimensional vector space.
Elements of order two in the Heisenberg group will define the spin 
representation of $so(2g+2)$. Restriction to suitable subspaces will give
the half spin representations and their relation with $G_g$ (a 
subquotient of $G_{g+1}$).
This relation between $so(2g+2)$ and the Heisenberg group is used in 
\ref{explop} to prove the Heisenberg invariance of the flat connections 
introduced in \ref{defcon}.

\subsubsection{} 
To accomodate both $G_g$ and $G_{g+1}$
we extend the construction of $\S$\ref{defB}. The inclusion 
$$
B:=\{1,\ldots,2g+2\}\hookrightarrow B^\sharp:=\{1,\ldots,2g+4\}
$$
and 
extension by zero of functions on $B$ to functions on $B^\sharp$ induces
$$
\tilde{F}_B:=
\left\{f:B\rightarrow\FF_2:\sum_{b\in B}f(b)=0\right\}\;\hookrightarrow \;
\tilde{F}_{B^\sharp}:=
\left\{f:B\rightarrow\FF_2:\sum_{b'\in B^\sharp}f(b')=0\right\}.
$$
The function $g'\in \tilde{F}_{B^\sharp}$ with $g'(b')=1$ 
(all $b'\in B^\sharp$) does not lie in
$\tilde{F}_B$, thus
$$
\tilde{F}_B
\hookrightarrow 
F_{B^\sharp}:=\tilde{F}_{B^\sharp}/\{0,g'\}.
$$
The image in $F_{B^\sharp}$ of the function $g
\;(\in \tilde{F}_{B})$ with $g(b)=1$ 
(all $b\in B$) is the element 
$$
p':=
x_{\{1,2,\ldots,2g+2\}}
=x_{\{2g+3,2g+4\}}\in F_{B^\sharp}.
$$ 
Using the definition of the symplectic form on $F_{B^\sharp}$ 
($\S$\ref{notiso}), which we denote by $E^\sharp$, one finds that:
$$
\tilde{F}_B=(p')^\perp,\qquad{\rm with}\quad 
p'^\perp:=\{x\in F_{B^\sharp}:\;E^\sharp(x,p')=0\;\}\quad
{\rm and}\quad F_B=(p')^\perp/\{0,p'\}.
$$
{}From $\S$\ref{notiso} we have an identification:
$$
F_{B^\sharp}=\FF^{g+1}\times\FF^{g+1},\qquad 
p'=(0,(1,1,\ldots,1))\qquad(\in F_{B^\sharp}).
$$
The Heisenberg group defined by this identification 
(cf.\ $\S$\ref{Hdefs}) will be 
denoted by $G^\sharp$, its Schr\"odinger representation by $U^\sharp$ (on the
$2^{g+1}$-dimensional vector space 
$V^\sharp:=\{f:\FF^{g+1}_2\rightarrow\CC\}$ 
with basis of $\delta$-functions $Y_{\sigma'}$, $\sigma'\in\FF_2^{g+1}$).

For any $x'\in (p')^\perp\;(\subset F_{B^\sharp})$,
the maps $U^\sharp_{p'}$ and $U^\sharp_{x'}$ 
commute. Therefore the $U^\sharp_{x'}$ with $x'\in (p')^\perp$ 
act on the two eigenspaces of $U^\sharp_{p'}$ 
which are:
$$
V^\sharp_+:=
\langle\ldots,Y_{\sigma'},\ldots\rangle_{\sigma_1'+\ldots+\sigma'_{g+1}=0},
\qquad
V^\sharp_-:=
\langle\ldots,Y_{\sigma'},\ldots\rangle_{\sigma_1'+\ldots+\sigma'_{g+1}=1}.
$$
A quotient map from $(p')^\perp\;(\subset 
(\FF_2^g\times\FF_2)\times (\FF_2^g\times\FF_2))$ to $F_B=\FF_2^g\times
\FF_2^g$
with kernel $\{0,p'\}$ is:
$$
(p')^\perp\longrightarrow F_B=(p')^\perp/\{0,p'\}\cong \FF_2^{2g},\qquad
x':=((a,a_{g+1}),(b,b_{g+1}))\longmapsto x:=(a,\bar{b})
$$
with $\bar{b}=(b_1+b_{g+1},\ldots,b_g+b_{g+1})$. 
We define an isomorphism of vector spaces:
$$
V^\sharp_+\longrightarrow V,\qquad
Y_{\sigma_1,\ldots,\sigma_g,\sigma_{g+1}}\longmapsto
X_{\sigma_1,\ldots,\sigma_g}.
$$

\begin{lem}
For any $x'\in (p')^\perp$ mapping to $x\in F_B$
there is a commutative diagram:
$$
\begin{array}{rcl}
V^\sharp_+&\stackrel{\cong}{\longrightarrow}&V\\
U^\sharp(t,x')\,\Big\downarrow&&\Big\downarrow\, U(t,x)\\
V^\sharp_+&\stackrel{\cong}{\longrightarrow}&V
\end{array}
$$
Moreover, if $1\leq j,k\leq 2g+1$,
then restriction of $U^\sharp_{x_{jk}}$ to $V^\sharp_+\cong V$ is given by
$\pm U_{x_{ij}}$.
\end{lem}

\begin{proof}
The two compositions along the square are (with $a,\,b,\,\sigma\in\FF_2^g$):
$$
Y_{\sigma,\sigma_{g+1}}\longmapsto
X_\sigma\longmapsto (t,(a,\bar{b}))X_\sigma=
t(-1)^{(a+\sigma)\bar{b}}X_{a+\sigma},\qquad{\rm and}
$$
$$
\begin{array}{rcl}
Y_{\sigma,\sigma_{g+1}}&\longmapsto&
(t,(a,a_{g+1}),(b,b_{g+1}))Y_{\sigma,\sigma_{g+1}}\\
&=&
t(-1)^{(a+\sigma)b+(a_{g+1}+\sigma_{g+1})b_{g+1}}
Y_{a+\sigma_g,a_{g+1}+\sigma_{g+1}}\\
&\longmapsto&
t(-1)^{(a+\sigma)b+(a_{g+1}+\sigma_{g+1})b_{g+1}}X_{a+\sigma}.
\end{array}
$$
Recall $E^\sharp(x',p')=0$, so
$a_{g+1}=\sum_{i=1}^ga_i$. As $Y_{\sigma,\sigma_{g+1}}\in V^\sharp_+$
we also have 
$\sigma_{g+1}=\sum_{i=1}^g\sigma_i$, thus:
$$
(a+\sigma)b+(a_{g+1}+\sigma_{g+1})b_{g+1}=
\left(\sum_{i=1}^g(a_i+\sigma_i)b_i\right)+
\left(\sum_{i=1}^g(a_i+\sigma_i)b_{g+1}\right)=(a+\sigma)\bar{b},
$$
which shows that the diagram commutes.
The last statement follows from the fact that such a $x'=x_{jk}$ lies in 
$(p')^\perp\;(\subset F_{B^\sharp})$ and that the homomorphism 
$x'\mapsto x\;(\in F_B)$ maps $x_{jk}$ to $x_{jk}$ (use \ref{notiso}).
Since $U^\sharp_{x_{jk}}=U^\sharp(t',x'),\;U_{x_{jk}}=U(t,x)$ 
for some $t',\,t\in\CC$ choosen such that each transformation
has order two, the statement 
follows from the commutativity of the diagram.
\end{proof}

\subsubsection{The Clifford algebra}
The Clifford algebra $C(Q)$ of the quadratic form 
$Q=x_1^2+\ldots+x_{2g+2}^2$ on
$V_Q=\CC^{2g+2}$ is the quotient of the tensor algebra $T(V_Q)$ by the 
two-sided ideal $I$ generated by the elements $v\otimes v-Q(v)$ for $v\in 
V_Q$ (\cite{FH}, $\S$20.1):
$$
C(Q):=T(V_Q)\Big/I=
(\CC\,\oplus\, V_Q\,\oplus\,V_Q\!\otimes\! V_Q\oplus\ldots)\Big/
(\ldots,-Q(v)+v\otimes v,\ldots)_{v\in V_Q}.
$$
Let $e_1,\ldots,e_{2g+2}$ be the standard basis vectors of $V_Q$, then 
the $C(Q)$ is generated by $\CC$ and the $e_j$ with relations:
$$
e_j\cdot e_j=1,\qquad e_j\cdot e_k+e_k\cdot e_j=0\quad(j\neq k),
$$
here $\cdot$ stands for the product induced by $\otimes$ on $C(Q)$.
The (associative) algebra $C(Q)$ becomes a Lie algebra by defining, as 
usual,
the  Lie bracket to be
$[x,y]:=x\cdot y-y\cdot x$. 

The following proposition relates the Heisenberg group, the Clifford 
algebra and the spin representation of the Lie algebra $so(2g+2)$.

\subsubsection{Proposition.} \label{halfspin}
With the notation as above we have:
\begin{enumerate}
\item
The $\CC$-linear map
$$
\gamma^\sharp:C(Q)\longrightarrow \mbox{End}(V^\sharp),\qquad
\lambda\mapsto\lambda I,\quad
e_k\mapsto U^\sharp_{x_{k,\,2g+4}}\qquad(\lambda\in\CC,\;\;k=1,\ldots,2g+2)
$$
defines an isomorphism of (associative $\CC$-) algebras.
\item
The linear map:
$$
\rho_s^\sharp:so(2g+2)\longrightarrow  
\mbox{End}(V^\sharp)\cong_{\gamma^\sharp} C(Q),\qquad
F_{jk}\longmapsto U^\sharp_{x_{j,2g+4}}U^\sharp_{x_{k,2g+4}}
\;(=\gamma^\sharp(e_j\cdot e_k))
$$
is an injective homomorphism of Lie algebras.
\item
The subspace $V^\sharp_+$ of $V^\sharp$ is invariant under the action of 
$so(2g+2)$. The Lie algebra representation
$$
\rho_s:so(2g+2)\longrightarrow \mbox{End}(V^\sharp_+)\cong 
\mbox{End}(V),\qquad
\rho_s(x):=\left.\rho_s^\sharp(x)_{\phantom{y}}\!\right|_{V^\sharp_+}
$$
is an (irreducible) half spin representation of $so(2g+2)$.
In particular, $V\cong V(\omega_{g+1})$.
\item
We have
$$
\rho_s(F_{jk})=\pm iU_{x_{jk}},\qquad(F_{jk}\in so(2g+2),\quad i^2=-1),
$$
where the $U_{x_{jk}}\;(\in \mbox{End}(V))$ with $x_{jk}\in \FF_2^{2g}$
are in the Schr\"odinger 
representation of the Heisenberg group $G$.
\end{enumerate}

\begin{proof}
By definition we have $(U^\sharp_{x'})^2=I$ for all $x'\in\FF_2^{2g+2}$. 
Moreover, 
$E^\sharp(x_{\{k,2g+4\}},x_{\{l,2g+4\}})=1$ when $k\neq l$ so
the corresponding maps anti-commute:
$U^{\sharp}_{x_{k,2g+4}}U^\sharp_{x_{l,2g+4}}=
-U^{\sharp}_{x_{l,2g+4}}U^{\sharp}_{x_{k,2g+4}}$. 
Then the map $\gamma^\sharp$ preserves the relations in $C(Q)$ and is thus an 
algebra homomorphism.
It is surjective because the matrices $U^\sharp_{x'}$, where $x'$ 
runs over $\FF_2^{g+1}\times\FF_2^{g+1}$ are a basis of ${\rm End}(V^\sharp)$ 
(the Schr\"odinger representation being irreducible) and any $U^\sharp_{x'}$ 
is a product (up to scalar multiple) of suitable $U^\sharp_{x_{j,2g+4}}$.
Therefore $\gamma^\sharp$ is 
an isomorphism since both algebras have the same dimension.

\noindent
(2)$\quad$
This is worked out in \cite{FH}, Lemma 20.7.
With their notations and our $Q$, we have 
$$
\phi:\wedge^2V_Q\stackrel{\cong}{\longrightarrow} so(2g+2),\qquad
e_j\wedge e_k\longmapsto F_{jk},
$$
$$
\psi:\wedge^2V_Q\longrightarrow C(Q),\qquad
e_j\wedge e_k\longmapsto e_j\cdot e_k\qquad(j\neq k),
\qquad{\rm and}\quad
\rho_s^\sharp=\psi\circ \phi^{-1}.
$$

\noindent
(3)$\quad$
By definition of $\rho_s^\sharp$ and the fact that $U^\sharp$ is a 
representation we have:
$$
\rho_s^\sharp(F_{jk})=U^\sharp_{x_{j,2g+4}}U^\sharp_{x_{k,2g+4}}=
c_{jk}U^\sharp_{x_{jk}},
$$
with 
$1\leq j,k\leq 2g+2$ and some $c_{jk}\in\CC^*$.
Thus the $\rho_s^\sharp(F_{jk})$ commute with $U^\sharp_{p'}$ in 
${\rm End}(V^\sharp)$. Therefore we obtain a Lie algebra representation on 
$V^\sharp_+$
cf.\ \cite{FH}, Prop.\ 20.15 (and identify their $End(\wedge^{even}W)$ 
with our $End(V^\sharp_+)$), where also the irreducibility is proved.

\noindent
(4)$\quad$We have $\rho_s^\sharp(F_{jk})=
U^\sharp_{x_{j,2g+4}}U^\sharp_{x_{k,2g+4}}$, thus 
$\rho_s^\sharp(F_{jk})^2=-I$  since
these two elements anti-commute. Therefore $\rho_s^\sharp(F_{jk})=\pm i 
U^\sharp_{x_{jk}}$, which acts as $\pm iU_{x_{jk}}$ on $V^\sharp_+$.
\end{proof}

\subsection{The Hitchin connection}\label{hitcon}

\subsubsection{} \label{PPdescent}
The symmetric group $S_{2g+2}$ acts on $\cP$ by permuting the coordinates,
the quotient wil be denoted by $\bP$.
The fundamental group of $\cP$ is the pure braid group and 
$\pi_1(\bP)=B_{2g+2}$, the Braid group.

In case $g=2$, we have $S_6\cong Sp(4,\FF_2)$, in fact there is a 
surjective map (which factors over the mapping class group)
$$
\pi_1(\overline{\cP})=B_{6}\longrightarrow Sp(4,\ZZ),
$$
the kernel of the composition $B_{6}\rightarrow Sp(4,\ZZ)\rightarrow
Sp(4,\FF_2)$ is $\pi_1(\cP)$.

>{}From the theory of theta functions we know that the groupscheme $\cG$
is trivialized on the cover $\tilde{\cP}$ of $\cP$ defined by the kernel 
of the composition:
$$
\phi:\pi_1(\overline{\cP})\longrightarrow
Sp(4,\ZZ)\longrightarrow Sp(4,\ZZ)/\Gamma_2(2,4)\cong A(G).
$$
Here we use Igusa's notation:
$$
\Gamma_2(2,4):=\left\{\pmatrix{I+2A&2B\cr 2C&I+2D}\in Sp(4,\ZZ):\;
diagonal(B)\equiv diagonal(C)\equiv (0,0)\;{\rm mod}\;2\right\}.
$$

Thus we consider the following diagram of etale Galois 
coverings:
$$
{\wP}\longrightarrow \cP\longrightarrow \overline{\cP}:=\cP/S_{6},
$$
and the corresponding exact sequence of covering groups:
$$
0\longrightarrow \FF_2^{2g}\longrightarrow A(G)\longrightarrow S_6
\longrightarrow 0.
$$
Recall that we have defined a projective representation 
$\tilde{T}:A(G)\rightarrow PGL(V)$ in \ref{projrep}.
This induces projective representations on each $S^kV$.

\begin{thm}\label{thmhitcon}
With the notation as in \ref{PPdescent}
 (except for $S_k$) and $g=2$ we have:
\begin{enumerate}
\item
There is a line bundle $\cN$ on ${\wP}$
such that the pull-back of $p_*\cL^{\otimes k}$ is isomorphic
to $S_k\otimes_\CC\cN$.
Here $S_k=S^kV(\omega_{g+1})$ 
and $V(\omega_{g+1})$ is a half spin representation of $so(6)\cong sl(4)$
(which is the standard representation (or its dual) of $sl(4)$).
\item\label{hits}
The Hitchin connection on the pull-back
of ${p}_*{\cal L}^{\otimes k}$ to $\wP$
is given by the pull-back to $\wP$ of
the connection on $S_k\otimes_\CC\cO_{\PP}$ defined by the one form
$$
{-1\over{16(k+2)}}\sum_{i,j\,i\neq j}
{{\sigma\tilde{\rho_s}(\Omega_{ij})\dd z_i}\over{z_i-z_j}}
$$
with $\tilde{\rho_s}:U(so(2g+2))\rightarrow{\rm End}(S_1)$ the half spin 
representation (then the $\sigma\tilde{\rho_s}(\Omega_{ij})$ give 
endomorphisms of each $S_k$, cf.\ \ref{symcon}).

\item \label{hitbp} 
The Hitchin connection on the bundle $S_k\otimes\cO_\wP$ over $\wP$
descends to a projective flat connection on the the bundle
$(S_k\otimes\cO_\wP)/A(G)$ over $\bP$, where $A(G)$ acts on $S_k$ via
the $k$-th symmetric power of its projective representation on $S_1$
and $A(G)$ acts on $\wP$ as $Gal(\wP/\bP)$.
This descent gives the
kth Hitchin connection of the natural descend of the family
of curves described in 1.4.
\end{enumerate}
\end{thm}

\begin{proof}
For the first part it suffices to show that the pull-back of 
$p:\PP\rightarrow\cP$ to $\tilde{\cP}$ is trivial.
Since ${\cal G}$ is isomorphic to the constant group scheme $G$ (cf.\ 
\ref{Hdefs} (with here $t\in\CC^*$)), the uniqueness of the Schr\"odinger 
representation gives the global triviallization.
We already observed that the line bundle $\cN$ does not interfere with the 
projective connections.

\noindent
(2)$\quad$ We use Theorem \ref{269} so we must verify $(Eq.)$ from \ref{crit}.
The covariant derivatives given are defined by the heat operator
$$
X=-
{1\over{16(k+2)}}\sum_{j\neq i}
{{\sigma\tilde{\rho_s}(\Omega_{ij})}}\otimes{{{\rm d}z_i}\over{z_i-z_j}}
\qquad(\in (S_2\otimes_\CC S_2^*)\otimes\Omega_\cP^1).
$$
The Heisenberg-invariance of this operator follows from the facts that
$\Omega_{ij}=F_{ij}\otimes F_{ij}$, that $\rho_s(F_{kj})=\pm iU_{x_{jk}}$
with $i^2=-1$
and that $U_yU_xU_y^{-1}=\pm U_x$.
Since $S^4(\tilde\rho_s)$ is a subrepresentation of 
$\tilde\rho:so(2g+2)\rightarrow {\rm End} S^2(\wedge^2\CC^{2g+2})$
we have 
$\sigma\tilde{\rho_s}(\Omega_{ij})=\sigma\tilde{\rho}(\Omega_{ij})$ on 
$S_4$.
For the connection $\nabla_0$ on $\tilde{p}_*\cO_{\PP}(1)= 
 S_1\otimes \cO_{\tilde{P}}$ 
we simply take $\nabla_0(fw)=w\otimes \dd f$. 
Then the equation $(Eq.)$ is exactly the statement of Theorem \ref{thmcon}.

\noindent
(3)$\quad$
We have a smooth projective family
of genus 2 curves over the scheme $\bar\cP$ given by
$$
     y^2 = f(x) = x^6+a_1 x^5 +a^2 x^4 + ... + a_6 = \prod_i (x-z_i),
$$
where $a_i$ is the $i$-th symmetric function of the $z_j$.
Thus by the general theory of Section 2 we have a ${\bf P}^3$-bundle
$\bar p : \bar {\bf P} \longrightarrow \bar\cP $
and a projective Hitchin connection $\bar{\bf D}_k$ on $\bar p_\ast \cO(k)$.
Let $\varphi : \tilde\cP \longrightarrow \bar\cP$ be
the natural map. Then there is an isomorphism
$$
(*)  \qquad\qquad     \varphi^* (  \bar {\bf P} ) = {\bf P}^3 \times \tilde\cP
$$
which is compatible with the action of $A(G)$ (which lies over the action
of $A(G)$ on $\tilde\cP$) on both sides: this isomorphism is given by the
trivialization of the action of the theta-group scheme $\cG$ which was
mentioned earlier. The isomorphism (*) induces an isomorphism
$$
(**) \qquad\qquad \varphi^* (  {\bar p}_*O(k) ) \cong S_k \otimes O_{\tilde\cP}
$$
which is compatible with the natural projective action of $A(G)$
on the LHS (via $\cG$) and the action of $A(G)$ on the RHS via the $Sym^k$
of its natural representation on $S_1$. The Heat operator given by $X$
(as in proof of \ref{thmhitcon}.\ref{hits}, see \ref{crit}) 
on the LHS of (*) is VIA (*) compatible
with
the pullback of the Hitchin heat operator of the LHS of (*). But both heat
operators
are invariant under the action of $A(G)$. Thus the projective connections
on the RHS and LHS of (**) agree and are compatible with the projective
action of $A(G)$ on both sides. In other words the natural descent
datum on the LHS of (**) agrees with the descent datum on the RHS
of (**).
\end{proof}

\subsubsection{Examples}\label{explop}
We give some examples of the $\rho_s(\Omega_{ij})$'s in case $g=2$. 
We identify $V=S_1$, with its standard basis $X_\sigma,\; 
\sigma\in\FF_2^2$. A linear map with matrix $(a_{kl})$ is then also given
by the linear differential operator $\sum_{kl}a_{kl}X_k\partial_l$ (since
$\partial_l(X_m)=0$ for $l\neq m$ and is $1$ if $l=m$).
Recall from 
\ref{notiso} that (for $g=2$):
$$
x_{12}=((0,0),(1,0)),\qquad{\rm so}\quad
U_{x_{12}}=X_{00}\partial_{00}+X_{01}\partial_{01}-
X_{10}\partial_{10}-X_{11}\partial_{11}.
$$
Similarly, we had $x_{16}=((1,0),(0,0))$ and $x_{26}=((1,0),(1,0))$ and so:
$$
U_{x_{16}}=X_{10}\partial_{00}+X_{11}\partial_{10}+
X_{00}\partial_{10}+X_{01}\partial_{10},\quad
U_{x_{26}}=i(-X_{10}\partial_{00}-X_{11}\partial_{10}+
X_{00}\partial_{10}+X_{01}\partial_{10}).
$$
Since $U_{x_{jk}}=\pm i\rho_s(F_{jk})$  and $U_{x_{jk}}^2=I$
we get (cf.\ \ref{symcon}):
$$
\sigma(\tilde{\rho_s}(\Omega_{12}))=
-2(X_{00}^2\partial_{00}^2+X_{01}^2\partial_{01}^2+
X_{10}^2\partial_{10}^2+X_{11}^2\partial_{11}^2+
2X_{00}X_{01}\partial_{00}\partial_{01}-2X_{10}X_{11}\partial_{10}\partial_{11}),
$$

\subsection{Local Monodromy}\label{localmon}

\subsubsection{Introduction}
We want to obtain some information on the representation of the mapping 
class group for genus two curves defined by Hitchin's connection. This 
representation has been studied, for arbitrary $g$ and $k$,
 by Moore and Seiberg \cite{MS},
Kohno \cite{K2} and in the case $k=2$ by Wright \cite{Wr}. 
In some sense, we merely find
(weaker, local) results which do agree with their results.

\subsubsection{The method} \label{monex}
We recall how to determine the local monodromy on the trival bundle
$\cO_\cP\otimes_\CC W$ with a connection $\nabla$.
We consider a holomorphic map
$$
\phi:D^\ast\longrightarrow \cP,\qquad{\rm with}\quad
D^*:=\{t\in\CC-\{0\}:\;|t|<\epsilon\,\}.
$$
(for small, positive $\epsilon$) and pull-back the
one form $\omega\in {\rm End}(W)\otimes\Omega^1_\cP$ 
which defines the connection (so $\nabla(fw)=w\dd f-f\omega(w)$. We write
$$
\phi^*\omega=
\left({R\over t}+A(t)\right){\rm d}t,\qquad R\in {\rm End}(W)
$$
with $A(t)$ holomorphic for $t=0$, and $R$ is called the residue at $t=0$.
We identify $\pi_1(D^*)=\ZZ$; $1\in\ZZ$ represents
a small circle traversed in anti-clockwise direction.
The eigenvalues of the monodromy of $1\in\pi_1(D^*)$ 
are the
$exp(-2\pi i\mu_j)$
were  the $\mu_i$ are the eigenvalues of $R$ on $W$.

In case the monodromy transformations are semi-simple, the local monodromy 
is conjugated to $exp(-2\pi iR)$ (cf. \cite{D}, p. 54).
However we could not prove the semi-simplicity
(but it seems to be known to the physicists).

To get the eigenvalues of the monodromy of the Hitchin connection
for a $\gamma \in\pi_1({\bP})$,
let $n$ be the order of its image in
$Gal(\cP/\bP)$ and let $\bar{\gamma}$ be its image in 
$A(G)=Gal({\wP}/{\bP})$. We choose a $\phi:D^*\rightarrow \cP$
in such a way that $\phi_*(1)=\gamma^n\;(\in \pi_1(\cP))$.
Then these eigenvalues are the eigenvalues of the matrix $exp(
\mbox{$-2\pi i\over n$}R)
\tilde{T}_{\bar{\gamma}}$, with $R$ determined as above for the
corresponding connection on $\cP$ (cf.\ \cite{K}). 
Since $\tilde{T}$ is a projective representation
(on $S_1$ and thus on any $S_k$), the set of eigenvalues are only defined up to
multiplication by one non-zero constant. This corresponds to the fact 
that we only have a projectively flat connection.

\subsubsection{Non-seperating vanishing cycle}\label{nsep}
Let $\gamma\in\pi_1(\bar{\cP})$ such that $\gamma^2=\phi_*(1)$ with
$$
\phi:D^\ast\longrightarrow \cP,\qquad
t\longmapsto (t+z_2,z_2,\ldots,z_6)
$$
where we fix distinct $z_j\in\CC$. (Then $\gamma$ corresponds to a Dehn 
twist in a simple non-seperating loop in the mapping class group.)
The residue of $\phi^*\omega$ for 
the connection on $S_k\otimes\cO_\cP$ corresponding to 
the Hitchin connection on $p_*\cL^{\otimes k}$ is
(cf.\ \ref{thmhitcon} and \ref{explop}):
$
\mbox{${-1\over {16(k+2)}}$}\sigma\tilde{\rho_s}(\Omega_{12}),
$
$$
-\tilde{\rho_s}(\Omega_{12})=
R_{12}:=
\left((X_{(0,0)}\partial_{(0,0)}+X_{(0,1)}\partial_{(0,1)})-
(X_{(1,0)}\partial_{(1,0)}+X_{(1,1)}\partial_{(1,1)})\right)^2
$$

\begin{lem}\label{spec12}
The eigenvalues of of the monodromy of $\gamma$ as in \ref{nsep}
on $S_k$ are:
$$
\lambda_c:=exp(-2\pi i \mbox{${c(c+1)}\over {k+2}$}),
\qquad
{\rm mult}(\lambda_c)=(k-2c+1)(2c+1)
$$
(up to multiplication by a non-zero constant independent of $c$),
with $2c\in\ZZ$ and $0\leq c\leq \mbox{$k\over 2$}$.
\end{lem}

\begin{proof}
Since the eigenvalues are only determined up to a constant and since
the difference between $2R_{12}$ and 
$\sigma(R_{12})$ is a multiple of the
Euler vector 
field, which acts by multiplication by $k$ on $S_k$, it suffices to consider 
the eigenvalues of $R_{12}$.
For a monomial in $S_k$ we have:
$$
R_{12}(X^{l_0}_{(0,0)}X^{l_1}_{(0,1)} 
X^{m_0}_{(1,0)}X^{m_1}_{(1,1)})=
=((l_0+l_1)-(m_0+m_1))^2(X^{l_0}_{(0,0)}X^{l_1}_{(0,1)} 
X^{m_0}_{(1,0)}X^{m_1}_{(1,1)})
$$
Thus each monomial is an eigenvector.
Let $b:=m_0+m_1$. 
Since the monomial has degree $k$, its eigenvalue is
$(k-2b)^2$. The multiplicity of the eigenvalue is $(k-b+1)(b+1)$
(the dimension of
the space of homogeneous polynomials in $2$ variables of degree $c$ is 
$c+1$).
Thus the eigenvalue $2(k-2b)^2$ of $2R_{12}$
has multiplicity $(k-b+1)(b+1)$.

The image $\bar{\gamma}$ of $\gamma$ in $A(G)$ is
 $\tilde{T}_{x_{12}}$ which acts on such a monomial by
(cf.\ \ref{extt}):
$$
\tilde{T}_{x_{12}}(X^{l_0}_{(0,0)}X^{l_1}_{(0,1)}
X^{m_0}_{(1,0)}X^{m_1}_{(1,1)})
=
i^b(X^{l_0}_{(0,0)}X^{l_1}_{(0,1)}
X^{m_0}_{(1,0)}X^{m_1}_{(1,1)}).
$$
The eigenvalues of $\gamma$ on $S_k$
are then, with $c:=b/2$:
$$
\begin{array}{rcl}
exp({-2\pi i\over 2}
\left(\mbox{$+1\over {16(k+2)}$}{2(k-2b)^2}+\mbox{$b\over 
4$}\right))
&=& exp(-2\pi i\left(\mbox{${k^2-8ck+16c^2+8ck+16c}\over 
{16(k+2)}$}\right))\\
&=&exp(-2\pi i\left(\mbox{${k^2}\over {16(k+2)}$}\right))\,
exp(-2\pi i\left(\mbox{${c(c+1)}\over {k+2}$}\right)).
\end{array}
$$
\end{proof}

\subsubsection{Seperating vanishing cycle}
Now we consider a $\gamma\in\pi_1(\bar{\cP})$ such that 
$\gamma=\phi_*(1)$ with
$$
\phi:D^*\longrightarrow \cP,\qquad
t\longmapsto (tz_1,tz_2,tz_3,z_4,z_5,z_6),
$$
with distinct nonzero $z_i$'s.

The residue in $t=0$ of the connection on $S_k\otimes\cP$ corresponding 
to Hitchin's connection is
$$
\mbox{$-1\over{16(k+2)}$}R_{123}, \qquad{\rm with}\quad
R_{123}:=2\sigma(\tilde{\rho_s}(\Omega_{12}+\Omega_{13}+\Omega_{23})).
$$

\begin{lem} \label{xq} Let $g=2$.
\begin{enumerate}
\item
There are constants $\lambda_k\in\QQ$ such that (in ${\rm End}(S_k)$):
$$
R_{123}=16QX_Q+\lambda_kI
\quad{\rm with}\quad
QX_Q=(X_{00}X_{10}-X_{01}X_{11})
(\partial_{00}\partial_{10}-\partial_{01}\partial_{11}).
$$
\item
We define a subspace of $S_k$ by:
$$
V_k:=Kernel(X_Q:=\partial_0\partial_1-\partial_2\partial_3:\, 
S_k\longrightarrow S_{k-2}).
$$
Then the vector space $S_k$ is a direct sum:
$$
S_k=V_k\oplus QV_{k-2}\oplus Q^2V_{k-4}\oplus\ldots
$$
Moreover, each subspace $Q^lV_{k-2l}$ is an eigenspace of $R_Q$
with eigenvalue:
$$
\lambda_l=
l(k-l+1)\qquad
{\rm and}\quad \dim V_{k-2l}=(k-2l+1)^2.
$$
\item
The eigenvalues of of the monodromy of $\gamma$ on $S_k$ are 
$$
\lambda_l:=exp(-2\pi i \mbox{${l(l+1)}\over 
{k+2}$})\qquad
mult(\lambda_l)=(k-2l+1)^2
$$
(up to multiplication by a constant independent of $l$), 
with $l\in\ZZ$ and $0\leq l\leq k/2$.
\end{enumerate}
\end{lem}

\begin{proof}
Recall $x_{12}$ corresponds
to $((0,0),(1,0))\in\FF^2_2\times\FF^2_2$ and 
$x_{23}=x_{\{2,3,4,5\}}+
x_{\{4,5\}}$ corresponds to $((1,1),(0,0))$
thus $x_{13}=x_{12}+x_{23}$ corresponds to 
$((1,1),(1,0))$.
Since $(1,1)(1,0)=1$, we have a
$+$ sign in front of $\Omega_{13}$ below.
Let  $\tilde{\rho_s}(\Omega_{00})$ 
be the square of the Euler vector field (which acts
as $k^2I$ on $S_k$).
Then:
$$
\begin{array}{rcr}
\tilde{\rho_s}(\Omega_{00})&=&+(X_{00}\partial_{00}+X_{01}\partial_{01}+
X_{10}\partial_{10}+X_{11}\partial_{11})^2,\\
\tilde{\rho_s}(\Omega_{12})&=&-(X_{00}\partial_{00}+X_{01}\partial_{01}-
X_{10}\partial_{10}-X_{11}\partial_{11})^2,\\
\tilde{\rho_s}(\Omega_{23})&=&-(X_{11}\partial_{00}+X_{10}\partial_{01}+
X_{01}\partial_{10}+X_{00}\partial_{11})^2,\\
\tilde{\rho_s}(\Omega_{13})&=&+(X_{11}\partial_{00}+X_{10}\partial_{01}-
X_{01}\partial_{10}-X_{00}\partial_{11})^2.
\end{array}
$$
The sum of these operators, minus the
 (degree one) parts which act as constants on each $S_k$, is $4QX_Q$.
Remembering the factor $2$ in $R_{123}$ and 
 in our symbols, we find the first statement.

\noindent
(2)$\quad$
This is well-known. Let $f\in V_{n}$, so $f$ is
homogeneous of degree $n$ and
$X_Q(f)=0$. Then an elementary
computation gives for $l\geq 1$ (with variables $X_i$ and 
$Q=X_0X_1-X_2X_3$):
$$
\begin{array}{cl}
&(\partial_0\partial_1-\partial_2\partial_3)(fQ^l)\\
=&(\partial_0 f)(\partial_1 Q^l)+(\partial_0 Q^l)(\partial_1 f)-
(\partial_2 f)(\partial_3 Q^l)-(\partial_2 Q^l)(\partial_3 f)+
f(\partial_0\partial_1-\partial_2\partial_3)(Q^l)
\\
=&lQ^{l-1}(X_0\partial_0+X_1\partial_1+X_2\partial_2+X_3\partial_3)(f)+
l(l+1)fQ^{l-1}
\\=& l(n+l+1)fQ^{l-1}.
\end{array}
$$
Writing $n=k-2l$ we get, for all integers $l$ with $0\leq 2l\leq k$:
$$
QX_Q(fQ^l)=l(k-l+1))fQ^l,\qquad f\in V_{k-2l}.
$$
Hence each $V_{k-2l}Q^l$ is an eigenspace for $R_Q$. By induction (the 
cases
$k=0,\,1$ being trivial) we may
assume that $S_{k-2}=\oplus S_{k-2-2l}Q^{l}$. Therefore $QS_{k-2}=
\oplus Q^{l+1}S_{k-2(l+1)}\subset S_k$ is a direct sum of eigenspaces of
$R_Q$
and none of the eigenvalues of $R_Q$ on this subspace is zero
($k+l(l+1)>0$ for $k\geq 2,\, l\geq 0$). Thus $\ker(R_Q)\cap
QS_{k-2}=\{0\}$ and $R_Q$ induces an isomorphism on
$QS_{k-2}$. Therefore $S_k=Ker (R_Q)\oplus
Im(R_Q)=\oplus Q^lS_{k-2l}$.
The dimension of $V_k$ is then $\dim S_k-\dim
S_{k-2}={{k+3}\choose 3}-{{k+1}\choose 3}=(k+1)^2$.

\noindent
(3)$\quad$
The image of $\gamma$ in $A(G)$ is 
trivial.
The eigenvalues of $R_{123}=8QX_Q$ 
on $S_k$ are $16l(k-l+1)$. Thus the eigenvalues of $\gamma$ are (see
\ref{monex}):
$$
\begin{array}{rcl}
exp(-2\pi i\mbox{${-1}\over{16(k+2)}$} 16l(k-l+1))
&{=}&
exp(-2\pi i\left(\mbox{${l^2-l-l(k+2)+2l}\over{k+2}$}\right))\\
&=&
exp(-2\pi i \mbox{${l(l+1)}\over{k+2}$}).
\end{array}
$$
\end{proof}

\subsubsection{Kohno's results} 
The monodromy of the Hitchin connection has been studied by Moore-Seiberg 
and Kohno (among others), we will relate the results above to
those contained in \cite{K2} for $g=2$.

Let $\gamma$ be a tri-valent graph with two edges, so there are three 
edges meeting
in each vertex.
For $k\in\ZZ_{\geq 0}$ we define a finite set $b_k(\gamma)$ of functions 
on the edges of $\gamma$ by:
$$
f:Edges(\gamma)\longrightarrow \{0,\mbox{$1\over 2$},1,
\ldots,\mbox{$k\over 2$}\}
\quad{\rm with}\quad
\left\{\begin{array}{r}
f(e_i)+f(e_j)+f(e_k)\leq k,\\
|f(e_i)-f(e_j)|\leq f(e_k)\leq f(e_i)+f(e_j),\\
f(e_i)+f(e_j)+f(e_k)\in\ZZ,
\end{array}
\right.
$$
for any three edges $e_i,\,e_j,\;e_k$ meeting in a vertex.
The Verlinde space
$V_k$ is the $\CC$-vector space with basis
$b_k(\gamma)$, it has the same dimension as $S_k$. 

The graph $\gamma$ corresponds to a pants decomposition of a genus two 
Riemann surface. Each vertex is a pant, homeomorphic to $\PP^1$ minus 3 
points,
the edges correspond to these points, which in turn correspond to 
`vanishing cycles' on the Riemann surface.
The mapping class group has a projective representation on $V_k$. 
The Dehn twist in a vanishing cycle corresponding to an edge 
$e$ of $\gamma$ acts, up to a scalar multiple, conjugated to the diagonal 
matrix with entries (cf.\ \cite{K2}, p.\ 217; p.\ 214, (2-1))
$$
exp(-2\pi i f(e)(f(e)+1)/(k+2)\qquad(f\in b_k(\gamma)).
$$
Using this recipe we find the same results as before (see Lemmas 
\ref{spec12}
and \ref{xq}):

\begin{lem}\label{tmon} With the notation as above:
\begin{enumerate} 
\item
The eigenvalues of the Dehn twist in a non-seperating cycle on $V_k$
are 
$exp(-2\pi i v(v+1)/(k+2))$ with multiplicity
$(k-2v+1)(2v+1)$ where $v\in \{0,\mbox{$1\over 2$},1,
\ldots,\mbox{$k\over 2$}\}$.
\item
The eigenvalues of a Dehn twist in a seperating cycle on $V_k$ 
are $exp(-2\pi il(l+1)/(k+2))$ with multiplicity
$(k-2l+1)^2$ where $l$ is an integer with $0\leq l\leq k/2$.
\end{enumerate}
\end{lem}

\begin{proof}
Consider the graph $\gamma$ with two vertices joined by three edges
$e_1,\,e_2,\,e_3$.
We first determine the set $b_k(\gamma)$. The edge $e_1$ corresponds to
a non-seperating cycle and we determine $f(e_1)$ for all $f\in 
b_k(\gamma)$, this will give the first result.

Writing $f_i:=f(e_i)$ for $f\in b_k(\gamma)$, the conditions are the same
in each of the two edges. Assume that we have $f_1\leq f_2\leq f_3$.
The conditions on the $f_i$ are then equivalent to the following 3 
conditions:
$$
f_1+f_2+f_3\leq k,\qquad f_3\leq f_1+f_2,\qquad f_1+f_2+f_3\in\ZZ,\qquad
f_i\in\{0,\mbox{$1\over 2$},1,
\ldots,\mbox{$k\over 2$}\}
$$
In case $f_1+f_2+f_3=l\;(\leq k)$ and $f_3>l/2$, the second condition 
implies
 $f_1+f_2>l/2$, a contradiction. Thus the set $b_k(\gamma)$ is a disjoint 
 union of sets $c_l(\gamma)$ for $l\in\ZZ$ and $0\leq l\leq k$:
 $$
 c_l(\gamma):=\{f\in b_k(\gamma):\;f_1+f_2+f_3=l,\quad f_i\leq l/2\;\},
 \qquad
 b_k(\gamma)=\stackrel{.}{\cup}_{0\leq l\leq k} c_l(\gamma).
 $$ 

Next we observe that if $f\in c_l(\gamma)$ then $f_1+f_2+f_3=l$ implies
$f_1+f_2=l-f_3\geq l/2$ (since $f_3\leq l/2$), thus certainly $f_3\leq 
f_1+f_2$, so the second condition is always fulfilled.

It remains to determine the number of triples $(f_1,f_2,f_3)$ with
$f_i\in \{0,\mbox{$1\over 2$},1,\ldots,\mbox{$l\over 2$}\}$
and $f_1+f_2+f_3=l$. Equivalently, we have to find the number of monomials 
$X^aY^bZ^c$ with $a+b+c=2l$ and $0\leq a,\,b,\,c\leq l$.
If such a monimial with $a+b+c=2l$ does not satifisfy $0\leq a,\,b,\,c\leq 
l$, then at most one of $a,\,b,\,c$ can be $>l$, say $c$, and then 
$a+b=2l-c$ which gives $(2l-c)+1$ possible couples $a,\,b$. Varying $c$ 
between $l+1$ and $2l$ we get $l+(l-1)+\ldots+1=(1/2)l(l+1)$ monomials 
with $c>l$ and thus there are $(3/2)l(l+1)$ monomials of degree $2l$ 
which do not satisfy the condition $0\leq a,\,b,\,c\leq l$. 
Then 
$$
|c_l(\gamma)|=\left({{2l+2}\atop 2}\right)-(3/2)l(l+1)={l+2\choose 2},
\quad{\rm hence}\quad
|b_k(\gamma)|=\sum_{l=0}^k {l+2\choose 2}=\left({{k+3}\atop 3}\right),
$$
so indeed $|b_k(\gamma)|=\dim S_k$.

Next we determine the number of triples $(v,f_2,f_3)\in 
b_k(\gamma)=\stackrel{.}{\cup}c_l(\gamma)$. 
First we will count such triples in $c_l(\gamma)$.
Let $a:=2v\in\ZZ_{\geq 0}$ with $a\leq l$, then
there are $2l-a+1$ couples $(b,c)$ with
$a+b+c=2l$. Of these, only $a+1$ have $b,\,c\leq l$ 
(they are $(l-a,l),\;(l-a+1,l-1),\ldots,(l,l-a)$).
This gives $a+1$ such triples in $c_l(\gamma)$ provided $f_1=a/2\leq l/2$.
Thus $v$ occurs for $l/2=v,\,v+1/2,\ldots,k/2$, that is for $k-2v+1$ 
values of $l/2$.
Thus there are  $(k-2v+1)(2v+1)$ triples $(v,f_2,f_3)$.

\noindent
(2)$\quad$
 We now consider the case of a seperating vanishing cycle.
 The recipe we folow is not explicitly given in \cite{K2}, but is similar
 to the one used above and gives a result that agrees with 
Lemma \ref{xq}.
 
Let $\delta$ be the graph with two vertices $v_1,\,v_2$
and three edges $e_1,\,e_2,\,e_3$ such that 
  begin and end of $e_i$ is $v_i$ and where $e_3$ connects $v_1$ and $v_2$.
Thus the edge $e_3$ corresponds to a seperating vanishing cycle.

As before, let $f_i:=f(e_i)$. The second set of inequalities on the $f(e_i)$ 
reduce to
$$
f_3\leq 2f_1,\qquad f_3\leq 2f_2;\qquad{\rm moreover}\quad
 f_1+f_1+f_3\in\ZZ\Rightarrow
f_3\in \ZZ\cap\{0,1/2,\ldots,k/2\},
$$
this is the third condition for the vertex $v_1$. Finally we have 
$2f_i+f_3\leq k$ so $f_3/2\leq f_i\leq (k-f_3)/2$ for $i=1,\,2$. 
Thus given $f_3$, we find $(k-2f_3)+1)^2$
possibilities for couples $f_1,\,f_2$ ($f_i\in\{
{{f_3}\over2},{{f_3+1}\over2},\ldots,
{{k-f_3}\over2}\}$). 

We observe that for fixed $f_1$, $f_3$ may have the values 
$0,1,\ldots,2f_1$,
correspondingly $f_2$ has $k+1,\; k+1-2,\ldots, k+1-4f_1$ different 
values.
Thus the multiplicity of $f_1$ is, again, 
$(2f_1+1)((k+1)-2(1+2+\ldots+2f_1)=(2f_1+1)(k+1-2f_1)$.
\end{proof}

\subsubsection{The cases $k=1,\;2$}
The case $k=1$ is particularly easy since the operators $\Omega_{ij}$ are 
homogeneous of degree two and thus act as zero on $S_1$. The monodromy 
representation then factors over the projective representation of 
$A(G)$. 

In case $k=2$, the space $S_k$ is a direct sum of 10 distinct, one dimensional 
$G$-representations (this is easy to verify, see also \cite{vG2}).
On each eigenspace, any $\sigma\tilde{\rho_s}(\Omega_{ij})$ acts as 
scalar multiplication by an integer. A (local) basis of flat sections 
of $S_2\otimes\cP$ is then given by 10 functions of the type
$(\prod(z_i-z_j)^{r_{ijm}})Q_m$ with $r_{ij}\in\QQ$ and $Q_m\in S_2$ is a  
basis of the eigenspace. 
These eigenspaces
are permuted by the projective representation of $A(G)$.

}

\vspace{\baselineskip}\noindent
Bert van Geemen,\\
Dipartimento di Matematica,\\
Universit\`{a} di Torino,\\
Via Carlo Alberto 10,\\
10123 Torino,\\
Italy.\\
geemen@dm.unito.it

\

\vspace{\baselineskip}\noindent
A.J.~de Jong,\\
Princeton University,\\
Department of Mathematics,\\
Fine Hall -- Washington Road,\\
Princeton, NJ 08544-1000,\\
USA.\\
dejong@math.Princeton.EDU

\end{document}